\documentclass[12pt]{article}
\usepackage{amscd,amsfonts,amssymb,amsmath,latexsym,array,hhline}
\usepackage[dvips]{graphics}
\makeatletter
\@addtoreset{equation}{section}
\renewcommand{\@begintheorem}[2]{\begin{trivlist}\it
\item[\hspace{\labelsep}{\bf #1\ #2.}]}
\renewcommand{\@opargbegintheorem}[3]{\begin{trivlist}\it
\item[\hspace{\labelsep}{\bf #1\ #2\ (#3).}]}
\renewcommand{\@endtheorem}{\end{trivlist}}
\renewcommand{\@cite}[2]{[{#1\if@tempswa ; #2\fi}]}
\makeatother

\newcommand{\paragr}{\hspace{6mm}}
\newtheorem{Theorem}{\paragr Theorem}[section]
\newtheorem{Lemma}[Theorem]{\paragr Lemma}

\newtheorem{Proposition}[Theorem]{\paragr Proposition}
\newtheorem{Definition}[Theorem]{\paragr Definition}

\newtheorem{Example}[Theorem]{\paragr Example}
\newtheorem{Corollary}[Theorem]{\paragr Corollary}
\newcommand{\Proof}{\texttt{Proof}. }
\newcommand{\Remark}{\texttt{Remark}. }

\newcommand{\Remarks}{\texttt{Remarks}. }

\newcommand{\End}{~\hfill $\Box\!$}

\newcommand{\skm}{\bigskip}
\newcommand{\skb}{\bigskip}

\newcommand{\hs}{\hspace*{2.5mm}}

\newcommand{\al}{\alpha}
\newcommand{\be}{\beta}
\newcommand{\ga}{\gamma}

\newcommand{\de}{\delta}
\newcommand{\De}{\Delta}

\newcommand{\la}{\lambda}
\newcommand{\si}{\sigma}

\newcommand{\eps}{\varepsilon}
\renewcommand{\phi}{\varphi}
\renewcommand{\kappa}{\varkappa}
\newcommand{\N}{\mathbb{N}}

\newcommand{\R}{\mathbb{R}}

\newcommand{\B}{\mathcal{B}}
\newcommand{\F}{\mathcal{F}}
\newcommand{\G}{\mathcal{G}}

\newcommand{\DDD}{\mathcal{D}}

\newcommand{\PPP}{\mathcal{P}}
\newcommand{\EX}{\mathcal{X}}

\newcommand{\EE}{\mathsf{E}}
\newcommand{\PP}{\mathsf{P}}
\newcommand{\QQ}{\mathsf{Q}}

\newcommand{\emp}{\emptyset}
\newcommand{\lb}{\langle}
\newcommand{\rb}{\rangle}
\newcommand{\wt}{\widetilde}
\newcommand{\wl}{\overline}
\newcommand{\xra}{\xrightarrow}
\newcommand{\da}{\downarrow}
\newcommand{\ua}{\uparrow}
\newcommand{\ds}{\displaystyle}
\newcommand{\cond}{\hspace{0.3mm}|\hspace{0.3mm}}

\newcommand{\Law}{\mathop{\rm Law}\nolimits}
\newcommand{\conv}{\mathop{\rm conv}\nolimits}

\newcommand{\pr}{\mathop{\rm pr}\nolimits}
\newcommand{\sgn}{\mathop{\rm sgn}\nolimits}

\newcommand{\essinf}{\mathop{\rm essinf}}

\newcommand{\argmin}{\mathop{\rm argmin}}

\newcommand{\Lin}{\mathop{\rm Lin}}
\newcommand{\cl}{\mathop{\rm cl}}
\newcommand{\cov}{\mathop{\sf cov}}
\newcommand{\var}{\mathop{\sf var}}
\newcommand{\corr}{\mathop{\sf corr}}
\renewcommand{\inf}{\mathop{\rm inf\rule[-0.8mm]{0mm}{1mm}}}

\newcommand{\Lea}{\Longleftrightarrow}

\newenvironment{mitemize}%
{\begin{list}{$\bullet$}{
\leftmargin=32pt
\rightmargin=0pt
\labelsep=5pt
\labelwidth=20pt
\itemindent=0pt
\topsep=5pt plus 2pt minus 4pt
\partopsep=2pt plus 1pt minus 1pt
\parsep=0pt
\itemsep=0pt}}%
{\end{list}}

\begin{document}
\begin{center}\bf
COHERENT MEASUREMENT OF FACTOR RISKS
\end{center}

\begin{center}\itshape\bfseries
Alexander S.~Cherny$^*$,\quad Dilip B.~Madan$^{**}$
\end{center}

\begin{center}
\textit{$^*$Moscow State University}\\
\textit{Faculty of Mechanics and Mathematics}\\
\textit{Department of Probability Theory}\\
\textit{119992 Moscow, Russia}\\
\texttt{E-mail: cherny@mech.math.msu.su}\\
\texttt{Webpage: http://mech.math.msu.su/\~{}cherny}
\end{center}

\begin{center}
\textit{$^{**}$Robert H.~Smith School of Business}\\
\textit{Van Munching Hall}\\
\textit{University of Maryland}\\
\textit{College Park, MD 20742}\\
\texttt{E-mail: dmadan@rhsmith.umd.edu}\\
\texttt{Webpage: http://www.rhsmith.umd.edu/faculty/dmadan}
\end{center}

\begin{abstract}
\textbf{Abstract.} We propose a new procedure for the
risk measurement of large portfolios.
It employs the following objects as the building blocks:
\begin{mitemize}
\item \textit{coherent risk measures}
introduced by Artzner, Delbaen, Eber, and Heath;
\item \textit{factor risk measures} introduced in this
paper, which assess the risks driven by particular
factors like the price of oil, S\&P500 index,
or the credit spread;
\item \textit{risk contributions} and
\textit{factor risk contributions}, which provide a
coherent alternative to the sensitivity coefficients.
\end{mitemize}
We also propose two particular classes of coherent risk
measures called \textit{Alpha V@R} and \textit{Beta V@R},
for which all the objects described above admit an
extremely simple empirical estimation procedure.
This procedure uses no model
assumptions on the structure of the price evolution.

Moreover, we consider the problem of the risk management
on a firm's level. It is shown that if the risk limits
are imposed on the risk contributions of the desks to
the overall risk of the firm (rather than on their
outstanding risks) and the desks are allowed to trade
these limits within a firm, then the desks automatically
find the globally optimal portfolio.

\medskip
\textbf{Key words and phrases.}
Alpha V@R,
Beta V@R,
capital allocation,
coherent risk measure,
extreme measure,
factor risk,
factor risk contribution,
law invariant risk measure,
portfolio optimization,
risk contribution,
risk management,
risk measurement,
risk sharing,
risk trading,
tail correlation,
Tail V@R,
Weighted V@R.
\end{abstract}

%=======================================================
\section*{Contents}

\noindent\hbox to \textwidth{\ref{I}\ \ Introduction
\dotfill \pageref{I}}

\noindent\hbox to \textwidth{\ref{CRM}\ \ Coherent Risk
\dotfill \pageref{CRM}}

\noindent\hbox to \textwidth{\ref{FR}\ \ Factor Risk
\dotfill \pageref{FR}}

\noindent\hbox to \textwidth{\ref{PO}\ \ Portfolio
Optimization \dotfill \pageref{PO}}

\noindent\hbox to \textwidth{\ref{RC}\ \ Risk Contribution
\dotfill \pageref{RC}}

\noindent\hbox to \textwidth{\ref{FRC}\ \ Factor Risk
Contribution \dotfill \pageref{FRC}}

\noindent\hbox to \textwidth{\ref{ORS}\ \ Optimal Risk
Sharing \dotfill \pageref{ORS}}

\noindent\hbox to \textwidth{\ref{SC}\ \ Summary and
Conclusion \dotfill \pageref{SC}}

\noindent\hbox to \textwidth{Appendix:
$L^0$-Version of the Kusuoka Theorem
\dotfill \pageref{A}}

\noindent\hbox to \textwidth{References
\dotfill \pageref{R}}

%========================================================
\section{Introduction}
\label{I}

\textbf{1. Overview.}
One of the basic problems of finance is:
\textit{How to measure risk in a proper way?}

Two most well-known and widely used in practice
approaches to this problem are variance and V@R.
However, both of them have serious drawbacks.
Variance penalizes high profits in the same way as
high losses.
Furthermore, the corresponding gain to loss ratio
$S(X)=\frac{\EE X}{\si(X)}$ known as the Sharpe ratio
does not have the monotonicity property:
$X\le Y$ does not imply that $S(X)\le S(Y)$.
In particular, an arbitrage possibility might have a
very low Sharpe ratio.
Concerning V@R, it takes into account only the
quantile of the distribution
without caring about what is happening to the left
and to the right of the quantile.
To put it another way, V@R is concerned only with the
probability of the loss and does not care about the size
of the loss.
However, it is obvious that the size of the loss should
be taken into account.
Let us remark in this connection that in the study of
the default risk the main two characteristics of a default
are its probability and its severity.
Further criticism of variance and V@R can be found
in~\cite{ADEH97} as well as in numerous discussions in
financial journals.

Recently, a new very promising approach to measure risk
was proposed in the landmark papers by Artzner, Delbaen,
Eber, and Heath~\cite{ADEH97}, \cite{ADEH99}
(these are the financial and the mathematical versions
of the same paper).
These authors introduced the notion of a
\textit{coherent risk measure}.
Since then, the theory of coherent risk measures has
rapidly been evolving; it already occupies a considerable
part of the modern financial mathematics.
Let us cite the papers~\cite{A02}, \cite{A04}, \cite{AT02},
\cite{ADEH99b}, \cite{C05e}, \cite{D02}, \cite{Dowd05},
\cite{FS02}, \cite{FS02b}, \cite{FRG02}, \cite{JMT04},
\cite{K01}, \cite{L04}, \cite{RU02}, \cite{T02},
to mention only a few.
Nice reviews on the theory of coherent risk measures
are given in~\cite{D05}, \cite[Ch.~4]{FS04}, \cite{S05}.
Much of the current research in this area deals with
defining properly a dynamic risk measure
(let us mention in this connection the
papers~\cite{CDK05}, \cite{DS05},
\cite{JR05}, \cite{R04}, \cite{RSE05}).
Currently, more and more research in the theory of
coherent risk measures is related to applications
to problems of finance rather than to the
study of ``pure'' risk measures.
In particular, the problem of capital allocation was
considered in~\cite{ADEH99b}, \cite{C061}, \cite{D05},
\cite{D01}, \cite{F03}, \cite{K05}, \cite{O99}, \cite{T02};
the problem of pricing and hedging was investigated
in~\cite{BE05'}, \cite{CGM01}, \cite{CGM02},
\cite{C061}, \cite{CK99}, \cite{D05},
\cite{JK01}, \cite{JR05}, \cite{LPST04},
\cite{N04}, \cite{RSE05}, \cite{Sek04}, \cite{S04};
the problem of the optimal portfolio choice was studied
in~\cite{C062}, \cite{RU00}, \cite{RUZ05};
the equilibrium problem was considered
in~\cite{BE05}, \cite{BE05'}, \cite{C062},
\cite{HK04}, \cite{JST05}, \cite{M04}.
This list is very far from being complete;
for example, on the Gloria Mundi web page over two
hundred papers are related to coherent risk measures.
The investigations mentioned above show that the whole
finance can be built based on coherent risk measures.
This is not surprising because risk
($\approx$~uncertainty) is at the very basis of finance,
and a new way of measuring risk yields new approaches
to all the basic problems of finance.
In some sources, the theory of coherent risk measures and
their applications is already called
the ``third revolution in finance'' (see~\cite{Sz04}).

\skb
\textbf{2. Coherent risk.}
A coherent risk measure is a functional defined on the
space of random variables that has the following form
\begin{equation}
\label{i1}
\rho(X)=-\inf_{\QQ\in\DDD}\EE_\QQ X.
\end{equation}
Here $X$ means the P\&L produced by some portfolio
over the unit time period (for example, one day)
and $\DDD$ is a set of probability measures
(all the measures from~$\DDD$ are assumed to be
absolutely continuous with respect to the original
measure~$\PP$).
From the financial point of view, $\rho(X)$ is the risk
of the portfolio.
Formula~\eqref{i1} has a straightforward economic interpretation:
we have a family~$\DDD$ of probabilistic scenarios;
under each scenario we calculate the average P\&L;
then we take the worst case.
Let us remark that typically a coherent risk measure is
defined as a functional on random variables satisfying
certain properties, and the representation theorem
states that it should have the form~\eqref{i1}.

The notion of a coherent risk measure is very convenient
theoretically, but when applying it to practice the
following problem arises immediately:
\textit{How to estimate $\rho(X)$ empirically?}
Of course, representation~\eqref{i1} cannot be used for
the practical calculation; furthermore, the class of
coherent risk measures ($=$ the class of different
sets~$\DDD$) is very wide. Thus, for the
practical purposes one needs to select a subclass of
coherent risk measures that admits an easy estimation
procedure. One of the best subclasses known so far is
\textit{Tail V@R}.
Tail V@R of order $\la\in(0,1]$ is the coherent risk
measure~$\rho_\la$ corresponding to
$\DDD=\{\QQ:d\QQ/d\PP\le\la^{-1}\}$.
It is easy to check that
$\rho_\la(X)=-\EE(X\cond X\le q_\la(X))$, where
$q_\la(X)$ is the $\la$-quantile of~$X$.
Thus, $\rho_\la$ is a very simple functional.
However, the sceptics can propose the following argument:
it is already hard to estimate $q_\la(X)$ due to rare tail
data and it is much harder to estimate the tail mean,
so the empirical estimation of~$\rho_\la$ is
problematic.

In this paper, we propose a one-parameter family of
coherent risk measures, which are extremely easy to
estimate. This is the class of risk measures of the form
$$
\rho_\al(X)=-\EE\min_{i=1,\dots,\al}X_i,
$$
where $\al$ is a fixed natural number and
$X_1,\dots,X_\al$ are independent copies of~$X$.
The parameter~$\al$ controls the risk aversion
of~$\rho_\al$: the higher is~$\al$, the more
risk averse is~$\rho_\al$.
We call $\rho_\al$ \textit{Alpha V@R}.
This family of risk measures is a very good substitute
for Tail V@R. It has the following advantages:
\begin{mitemize}
\item $\rho_\al$ is very intuitive;
\item $\rho_\al$ depends on the whole distribution
of~$X$ and not just on the tail as $\rho_\la$;
\item it is very easy to estimate $\rho_\al(X)$ from
the time series for~$X$.
\end{mitemize}
Let us remark that for large data sets
the empirical estimation of
$\rho_\al$ from time series is faster than
the empirical estimation of
$\text{V@R}_\la$ (!)
because it does not require the ordering of the time
series.
We believe that $\rho_\al$ is the best one-parameter
family of coherent risk measures.

Alpha V@R is a subclass of the two-parameter family of
coherent risk measures, which we call \textit{Beta V@R}.
This is the class of risk measures of the form
$$
\rho_{\al,\be}(X)=
-\EE\Bigl[\sum_{i=1}^\be X_{(i)}\Bigr],
$$
where $\be<\al$ are fixed natural numbers and
$X_{(1)},\dots,X_{(\al)}$ are the order statistics
obtained from the sequence $X_1,\dots,X_\al$ of
independent copies of~$X$.
In particular, $\rho_{\al,1}=\rho_\al$.
We believe that $\rho_{\al,\be}$ is the best
two-parameter family of coherent risk measures.

\skb
\textbf{3. Factor risk.}
The theory of coherent risk measures is already rather rich.
However, so far, an important issue has been absent:
\textit{How to measure the separate risks of a
portfolio induced by factors like the price of oil,
S\&P 500 index, or the credit spread?}
Measuring these risks separately is very important.
Suppose, for example, that a significant change in the
price of oil is expected in the near future.
Then having a big exposure to the price of oil
is dangerous.

Of course, one might argue at this point that $\rho(X)$
depends on $\Law X$, and this distribution already
takes into account high risk induced by possible oil
price moves (as well as all the other risks),
so that there is no need to consider separate risks.
Our reply to this possible criticism is as follows.
Suppose that we are trying to assess empirically
the risk of a large portfolio.
As different assets in a large portfolio have
different durations (like options or bonds),
a joint time series for all of them
does not exist.
So, for the empirical estimation we should consider
the main factors driving the risk, express the values
of all the assets through these factors, and use time
series for the factors.
Another advantage of factor risks is as follows.
When choosing data for the empirical risk estimation,
there is always a conflict between accuracy and flexibility
to the recent changes (the more data we take, the more
accurate are estimates, but less is the flexibility to
recent changes).
However, if we are looking for the separate risks
driven by separate factors, then this conflict can be
resolved. Namely, for a factor one can take time series
whose time step is stretched according to the current
volatility of the factor (for details, see
Section~\ref{FR}). This time change procedure enables
one to use arbitrarily large data sets and at the same
time immediately react to the volatility changes.

We propose a way for the coherent measurement of factor
risks. Let $X$ be the P\&L of a portfolio produced over a unit
time period and $Y$ be the increment of some factor/factors
over this period ($Y$ might be multidimensional).
The \textit{factor risk} of $X$ induced by $Y$ is
defined as
$$
\rho^f(X;Y)=-\inf_{\QQ\in\EE(\DDD\cond Y)}\EE_\QQ X,
$$
where $\DDD$ is the set standing in~\eqref{i1} and
$$
\EE(\DDD\cond Y):=\{\EE(Z\cond Y):Z\in\DDD\}.
$$
Here we identify measures from $\DDD$ with their
densities with respect to~$\PP$.
Thus, $\rho^f(\,\cdot\,;Y)$ is again a coherent risk measure.

It is easy to check that
$$
\rho^f(X;Y)=\rho(\EE(X\cond Y)).
$$
This expression clarifies the essence of factor risk:
$\rho^f(X;Y)$ takes into account the risk of~$X$ driven
only by~$Y$ and cuts off all the other risks.
The value $\rho^f(X;Y)$ might be looked at as the
coherent counterpart of the sensitivity of~$X$ to~$Y$.
However, an important difference from the sensitivity
analysis is the non-linearity of $\rho^f(X;Y)$.

The notion of factor risk can be employed in two forms:
\begin{mitemize}
\item[\bf 1.] We take one factor (i.e. $Y$ is
one-dimensional), and then $\rho^f$ measures the risk
induced by this factor.
\item[\bf 2.] We take all the main factors driving the
risk (i.e. $Y$ is multidimensional), and then $\rho^f$
measures the ``non-diversifiable'' risk and serves
as a good approximation to~$\rho$.
In other words, we are taking the main factors, express
the values of all the assets in the portfolio through
these factors, and look at the coherent risk.
This is a standard procedure of the risk measurement,
the only difference is that coherent risk measures are
employed.
\end{mitemize}

A very pleasant feature of $\rho^f$ is that it is easily
calculated for large portfolios.
Suppose, for example, that our portfolio consists of
$N$ assets, so that $X=X^1+\dots+X^N$.
Then
$$
\rho^f\Bigl(\sum_{n=1}^N X^n;Y\Bigr)
=\rho\Bigl(\sum_{n=1}^N f^n(Y)\Bigr),
$$
where $f^n(y)=\EE(X^n\cond Y=y)$.
In order to calculate this value, we need not know
the joint distribution of~$X^1,\dots,X^N$
(if the portfolio consists of several thousands assets,
finding this distribution is not a very pleasant problem).
Instead, we only need the joint distributions of
$(X^n,Y)$, $n=1,\dots,N$.
In essence, there is no difference whether our portfolio
consists of 1000 or 1000000 assets because different
$f^n$s are joined simply by the summation procedure.

For coherent factor risks based on Alpha V@R and Beta V@R,
we consider an empirical estimation procedure similar
to the historic V@R estimation (for the description of
these methods, see~\cite[Sect.~6]{M02}).
This procedure has the following advantages:
\begin{mitemize}
\item arbitrarily large data sets for factors are available;
\item for one-factor risks we can use the time change
procedure, which enables one to use arbitrarily large
data sets and immediately react to the volatility
changes;
\item the procedure is simple and works at the same speed
as the historic V@R estimation;
\item the use of coherent risk measures is much wiser
than the use of V@R;
\item the procedure is completely non-linear;
\item the procedure employs no model assumptions
(except, of course, for those used in the calculation
of $\EE(X\cond Y=y)$).
\end{mitemize}

\skb
\textbf{4. Portfolio optimization.}
The main message of the first part of the paper is:
it is reasonable to assess the risk of a position not
just as a number, but as a vector of factor risks
$\rho^f(\,\cdot\,;Y^1),\dots,\rho^f(\,\cdot\,;Y^M)$.
We also consider the problem of
portfolio optimization when the constraints on the
portfolio are given as $\rho^f(X;Y^m)\le c^m$,
$m=1,\dots,M$.
It turns out that this problem admits a simple geometric
solution that is similar to the one given
in~\cite[Subsect.~2.2]{C062}
for the case of a single constraint $\rho(X)\le c$.
We do not insist that this solution is the one that
should be implemented in practice, but it gives a nice
theoretic insight into the form of the optimal portfolio.
A possible practical approach to this problem is also
discussed.

\skb
\textbf{5. Risk contribution.}
The functional $\rho$ measures
the outstanding risk of a portfolio.
However, if a big firm assesses the risk of a trade,
it should take into account not the outstanding risk
of the trade, but rather its impact on the risk of the
whole firm. Thus, if the P\&L of the trade is~$X$
and the P\&L produced over the same period by the whole
firm is~$W$, the quantity of interest is $\rho(W+X)-\rho(W)$.
If $X$ is small as compared to~$W$, then a good
approximation to this difference is the
\textit{risk contribution} of~$X$ to~$W$.
The notion of risk contribution based on coherent risk
measures was considered in~\cite{D01}, \cite{F03},
\cite{K05}, \cite{O99}, \cite{T02}.
In this paper, we are taking the definition
proposed in~\cite[Subsect.~2.5]{C061}.
One of equivalent definitions of the coherent risk
contribution $\rho^f(X;W)$ of~$X$ to~$W$ is
$$
\rho^c(X;W)
=\lim_{\eps\da0}\eps^{-1}(\rho(W+\eps X)-\rho(W)).
$$
Note that if $X$ is small as compared to~$W$, then
$$
\rho(W+X)-\rho(W)\approx\rho^c(X;W).
$$
In most typical situations $\rho^c$ admits the
representation
$$
\rho^c(X;W)=-\EE_\QQ X,
$$
where $\QQ$ is the \textit{extreme measure}
defined as $\argmin_{\QQ\in\DDD}\EE_\QQ W$
($\DDD$ is the set standing in~\eqref{i1}).
Note that in this case $\rho^f$ is
linear in~$X$.

For the most important classes of coherent risk measures,
risk contribution admits a simple empirical estimation
procedure.
As shown in the paper, for Alpha V@R, $\rho^c$ has the
following form:
$$
\rho_\al^c(X;W)=-\EE X_{\argmin\limits_{i=1,\dots,\al}W_i},
$$
where $(X_1,W_1),\dots,(X_\al,W_\al)$ are independent
copies of $(X,W)$.
In particular, it is not hard to compute this quantity for
the typical case, where the firm's portfolio consists
of many assets, i.e. $W=\sum_{n=1}^N W^n$ with a large
number~$N$.
A similar simple representation is provided for the
Beta V@R risk contribution.

We also provide formulas for the empirical estimation
of risk contributions for the class of risk measures,
which we call \textit{Weighted V@R} (the term
\textit{spectral risk measures} is also used in the
literature); this is a wide class containing,
in particular, Beta V@R.

We also study the properties of the coefficient
$$
\kappa(X;W)=\frac{u^c(X;W)}{u(X)}.
$$
which measures the tail correlation between $X$ and $W$.

\skb
\textbf{6. Factor risk contribution.}
In this paper, we introduce the notion of
\textit{factor risk contribution}. It is simply the risk
contribution applied to the risk measure $\rho^f$, i.e
$$
\rho^{fc}(X;Y;W)
=\lim_{\eps\da0}\eps^{-1}(\rho^f(W+\eps X;Y)-\rho^f(W;Y)).
$$
The quantity $\rho^{fc}(\,\cdot\,;W;Y)$ may be viewed as
a coherent alternative to the sensitivity coefficient.
There are, however, two important differences:
\begin{mitemize}
\item The sensitivity coefficient is a nice estimate of
risk provided that the change in the corresponding factor
is small (i.e. $Y$ is small).
On the other hand, the factor risk contribution is a
nice estimate of risk provided that $X$ is small as
compared to~$W$ and there are no requirements on the
size of~$Y$.
\item The sensitivity coefficient measures the
sensitivity of a portfolio to a market factor.
On the other hand, the factor risk contribution measures
the sensitivity of a portfolio~$X$ both to the market
factor~$Y$ and to the firm's portfolio~$W$.
So, $\rho^{fc}(\,\cdot\,;Y;W)$ is a ``firm-specific''
coherent sensitivity to the factor~$Y$.
\end{mitemize}

As shown in the paper,
$$
\rho^{fc}(X;Y;W)=\rho^c(\EE(X\cond Y);\EE(W\cond Y)).
$$
This formula reduces the computation of $\rho^{fc}$ to
two steps: calculating the functions
$f(y)=\EE(X\cond Y=y)$, $g(y)=\EE(W\cond Y=y)$
and computing~$\rho^c$.

\skb
\textbf{7. Risk sharing.}
The relevance of factor risk contributions becomes clear
from our considerations of the risk sharing problem.
One of the basic problems of the central management of
a firm consisting of several desks is:
How to impose the limits on the risk of each desk?
Typically, this problem is approached as follows.
By looking at the performance of each desk, the central
management decides which desks should grow and which ones
should shrink and chooses the risk limits accordingly.
Nowadays, the procedure of choosing the risk limits is
to a large extent a political one.
At the same time, the central managers would like to
have a quantitative rather than political
procedure of choosing these limits.
An idea proposed by practitioners is:
Instead of giving each desk a fixed risk limit, it
might be reasonable to allow the desks to trade the
risk limits between themselves. For example, if one
desk is not using its risk limit completely, it might
sell the excess risk limit to another desk, which needs it.

In this paper we address the following problem:
\textit{Is it possible to arrange a
market of risk within a firm in such a way that the
resulting competitive optimum would coincide with the
global one?}
(By the global optimum we mean the one attained if the
central management had possessed all the information
available to all the desks and had been able to solve
the corresponding global optimization problem.)
The hope of the positive answer is justified by a
well-known equivalence established in the expected utility
framework between the global optimum (known also as the
Pareto-type or the soviet-type optimum) and the
competitive optimum (known also as the
Arrow--Debreu-type or the western-type optimum);
see~\cite{K90}.
A coherent risk counterpart of this result was
established in~\cite[Sect.~4]{C062}. The difference
between these results and our setting is that here the
objects traded are risk limits rather than financial
contracts. The results of this paper (they are
established within the coherent risk framework)
are as follows:
\begin{mitemize}
\item If the desks are measuring their outstanding risks
and are keeping them
within the risk limits, then the competitive optimum
is \textit{not} the global one.
\item If the desks are keeping track of their
risk contributions to the whole firm,
then the competitive optimum coincides with the global
one.
\end{mitemize}
Moreover, it turns out that the global optimum is achieved
regardless of what the initial allocation of risk
limits between the desks is.
(By the initial allocation we mean the risk limits
given to the desks by the central management before the
desks start to trade their risk limits.)
Our result applies not only to the case of one risk
constraint on the firm's portfolio,
but also to the case of several risk constraints
(a typical example is the constraints on each factor risk).

\skb
\textbf{8. Structure of the paper.}
In essence, the paper consists of two parts.
The first part (Sections~\ref{CRM}--\ref{PO})
deals with factor risk.
In Section~\ref{CRM}, we recall basic facts and examples
related to coherent risk measures.
Section~\ref{FR} deals with factor risk.
In Section~\ref{PO}, we study the problem of portfolio
optimization under limits imposed on factor risks.
The second part (Sections~\ref{RC}--\ref{ORS})
deals with (factor) risk contribution.
In Section~\ref{RC}, we recall basic facts
related to risk contribution.
Section~\ref{FRC} deals with factor risk contribution.
In Section~\ref{ORS}, we study the problem of risk
sharing between the desks of a firm.
The Appendix contains the $L^0$-version
of the Kusuoka theorem, which is needed for some
statements of the paper and is of interest by itself.
The recipes for the practical risk measurement are
gathered in Section~\ref{SC}, where we also compare
various empirical risk estimation techniques considered
in this paper with the classical ones like parametric V@R,
Monte Carlo V@R, and historic V@R.
The reader interested in practical applications only
might proceed directly to Section~\ref{SC}.

\bigskip
{\itshape\bfseries Acknowledgement.}
A.~Cherny expresses his thanks to D.~Heath and S.~Hilden for
the fruitful discussions related to the risk sharing problem.

%=======================================================
\section{Coherent Risk}
\label{CRM}

\textbf{1. Basic definitions and facts.}
Let $(\Omega,\F,\PP)$ be a probability space.
It is convenient to consider instead of coherent risk
measures their opposites called coherent utility
functions. This enables one to get rid of numerous
minus signs.
Recall that $L^\infty$ is the space of bounded random
variables on $(\Omega,\F,\PP)$.

The following definition was introduced in~\cite{ADEH97},
\cite{ADEH99}, \cite{D02}.

\begin{Definition}\rm
\label{CRM1}
A \textit{coherent utility function on $L^\infty$} is a
map $u\colon L^\infty\to\R$ satisfying the properties:
\begin{mitemize}
\item[(a)] (Superadditivity) $u(X+Y)\ge u(X)+u(Y)$;
\item[(b)] (Monotonicity) If $X\le Y$, then $u(X)\le u(Y)$;
\item[(c)] (Positive homogeneity) $u(\la X)=\la u(X)$ for
$\la\in\R_+$;
\item[(d)] (Translation invariance) $u(X+m)=u(X)+m$ for $m\in\R$;
\item[(e)] (Fatou property) If $|X_n|\le1$,
$X_n\xra{\PP}X$, then $u(X)\ge\limsup_n u(X_n)$.
\end{mitemize}
The corresponding \textit{coherent risk measure} is
$\rho(X)=-u(X)$.
\end{Definition}

\Remarks
(i) From the financial point of view, $X$ means the
P\&L produced by some portfolio over the unit time
period (taken as the basis for risk measurement)
and discounted to the initial time.
Actually, all the financial quantities in this paper
are the discounted ones.
However, the unit time period is typically small
(for example, one day), and for such
time horizons the discounted values are very close to the
actual ones. For this reason, below we skip the word
``discounted''.

(ii) The superadditivity property of~$u$
($=$~the subadditivity of~$\rho$) has the following
financial meaning: if we have a
portfolio consisting of several subportfolios and the
risk of each subportfolio is small, then the risk of
the whole portfolio is small.
V@R satisfies all the conditions of the above definition
except for the subadditivity, and this leads to serious
drawbacks of V@R. This can be illustrated by a simple
example proposed in~\cite{ADEH97}.
Suppose that a portfolio consists of 25 subportfolios
(corresponding to several agents),
i.e. its P\&L $X$ equals $\sum_{n=1}^{25}X^n$.
Suppose that $X^n=I_{(A^n)^c}-100 I_{A^n}$, where
$\PP(A_i)=1/25$ and $A^1,\dots,A^{25}$ are disjoint sets
($(A^n)^c$ denotes the complement of~$A^n$).
This means that each agent employs a spiking strategy.
If risk is measured by $\text{V@R}_{0.05}$, then
the risk of each subportfolio is negative (meaning that
each subportfolio is extremely good from the viewpoint
of this risk measure), while the P\&L produced by the
whole portfolio is identically equal to $-76$!

A natural example, in which V@R is subadditive, is the
Gaussian case: if $X,Y$ have a jointly Gaussian distribution, then
$\text{V@R}_\la(X+Y)\le\text{V@R}_\la(X)+\text{V@R}_\la(Y)$.
But actually on the set of (centered) Gaussian variables
all the reasonable risk measures like V@R, variance,
any law invariant coherent risk measure (see Example~\ref{CRM7})
coincide up to multiplication by a positive constant.
However, for general random variables this is not the
case, V@R and variance exhibit serious drawbacks
(see~\cite{ADEH97}), and coherent risk measures are
really needed.\End

\skm
The theorem below is the basic representation theorem.
It was established in~\cite{ADEH99} for
the case of a finite $\Omega$ (in this case the axiom~(e)
is not needed) and in~\cite{D02} for the general case
(the proof can also be found
in~\cite[Cor.~4.35]{FS04} or~\cite[Cor.~1.17]{S05}).
We denote by $\PPP$ the set of probability measures
on~$\F$ that are absolutely continuous with respect to~$\PP$.
Throughout the paper, we identify measures from~$\PPP$
(these are typically denoted by~$\QQ$)
with their densities with respect to~$\PP$
(these are typically denoted by~$Z$).

\begin{Theorem}
\label{CRM2}
A function~$u$ satisfies
conditions {\rm(a)--(e)} if and only if
there exists a non-empty set $\DDD\subseteq\PPP$ such that
\begin{equation}
\label{crm1}
u(X)=\inf_{\QQ\in\DDD}\EE_\QQ X,\quad X\in L^\infty.
\end{equation}
\end{Theorem}

\Remarks
(i) Let us emphasize that a coherent risk measure is defined
on random variables and not on their distributions.
Using representation~\eqref{crm1}, it is easy to construct
an example of a risk measure $\rho$ and two random variables
$X$ and $Y$ with $\Law X=\Law Y$, but with
$\rho(X)\ne\rho(Y)$.
A particularly important subclass of coherent risk
measures is the class of \textit{law invariant} ones,
i.e. the risk measures $\rho$ that depend only on the
distribution of~$X$.
However, it would not be a nice idea to include law
invariance as the sixth axiom in Definition~\ref{CRM1}.
Indeed, the basic risk measure used by an agent is
typically law invariant.
But there are many ``derivative'' risk measures like
$\rho^f(\,\cdot\,;Y)$
(many examples of naturally arising ``derivative'' risk
measures can be found in~\cite{C061}, \cite{C062},
\cite{D05}), and these ones need not be law invariant
even if the basic risk measure~$\rho$ is law invariant.

(ii) Coherent risk measures are primarily aimed at measuring
risk. But they can be used to measure
the risk-adjusted performance as well.
The risk-adjusted performance based on coherent risk
is a functional of the form $p(X)=\EE X-\la\rho(X)$,
where $\rho$ is a coherent risk measure and $\la\in\R_+$.
(By $\EE$ we will always denote the expectation with
respect to the original measure~$\PP$.)
Note that the functional $X\mapsto\EE X$ is a coherent
utility (it is sufficient to take $\DDD=\{\PP\}$
in~\eqref{crm1}).
Furthermore, if
$$
u_n(X)=\inf_{\QQ\in\DDD_n}\EE_\QQ X,\quad n=1,2
$$
are two coherent utilities and $\la\in[0,1]$, then
$$
\la u_1(X)+(1-\la)u_2(X)
=\inf_{\QQ\in\la\DDD_1+(1-\la)\DDD_2}\EE_\QQ X
$$
is also a coherent utility
(we use the notation $\la\DDD_1+(1-\la)\DDD_2
=\{\la\QQ_1+(1-\la)\QQ_2:\QQ_n\in\DDD_n\}$).
Thus, $(1+\la)^{-1}p(X)$ is again a coherent utility,
which is a very convenient ``stability'' feature.
As a result, coherent utility/risk can be used
\begin{mitemize}
\item[\bf 1.] to measure risk;
\item[\bf 2.] to measure the risk-adjusted performance, i.e. utility.
\end{mitemize}

(iii) Coherent utility may serve as a substitute for the
classical expected utility. The techniques like
utility-indifference pricing, utility-based optimization,
utility-based equilibrium, etc. can be transferred
from the expected utility framework to the coherent
utility framework.
Note that the intersection of these two classes of utility
is trivial: it consists only of the functional
$X\mapsto\EE X$ (note that all the other expected utilities
do not satisfy the translation invariance property).

(iv) A generalization of coherent utility is the notion
of \textit{concave utility} introduced by F\"ollmer and
Schied~\cite{FS02}, \cite{FS02b}. A functional
$u\colon L^\infty\to\R$ is a
\textit{concave utility function}
if it satisfies the axioms (b), (d), (e) of
Definition~\ref{CRM1} as well as the condition
\begin{mitemize}
\item[(a')] (Concavity)
$u(\la X+(1-\la)Y)\ge\la u(X)+(1-\la)u(Y)$ for
$\la\in[0,1]$.
\end{mitemize}
The corresponding \textit{convex risk measure} is
$\rho(X)=-u(X)$.
The representation theorem states that~$u$ is a concave
utility function if and only there exists a non-empty set
$\DDD\subseteq\PPP$ and a function $\al\colon\DDD\to\R$
such that
$$
u(X)=\inf_{\QQ\in\DDD}(\EE_\QQ X+\al(\QQ)),\quad X\in L^\infty
$$
(the proof can be found in~\cite{FS02},
\cite[Th.~4.31]{FS04}, or~\cite[Th.~1.13]{S05}).
A natural example of a concave utility function is
$$
u(X)=\sup\Bigl\{m\in\R:\inf_{\QQ\in\DDD}
\EE_\QQ U(X-m)\ge x_0\Bigr\},\quad X\in L^\infty,
$$
where $U:\R\to\R$ is a concave increasing function,
$\DDD\subseteq\PPP$, and $x_0\in\R$ is a fixed threshold.
This object is closely connected with the robust version
of the Savage theory developed by Gilboa and
Schmeidler~\cite{GS89}. For a detailed study of these
concave utilities, see~\cite{FS02b},
\cite[Sect.~4.9]{FS04}, or~\cite[Sect.~1.6]{S05}.

The theory of convex risk measures is now a rather large
field. However, in applications to problems of finance
like pricing and optimization, coherent risk measures
turn out to be much more convenient than the convex ones.
For this reason, we will consider only coherent risk
measures.\End

\skm
So far, a coherent risk measure has been defined on bounded
random variables.
Let us ask ourselves the following question:
Are ``financial'' random variables like the increment
of a price of some asset indeed bounded?
The right way to address this question is to
split it into two parts:
\begin{mitemize}
\item Are ``financial'' random variables bounded in practice?
\item Are ``financial'' random variables bounded in theory?
\end{mitemize}
The answer to the first question is positive
(clearly, everything is bounded by the number of the
atoms in the universe).
The answer to the second question is negative because
most distributions used in theory (like the lognormal
one) are unbounded.
So, as we are dealing with theory, we need to extend
coherent risk measures to the space $L^0$ of all random
variables.

It is hopeless to axiomatize the notion of a risk measure
on $L^0$ and then to obtain the corresponding
representation theorem.
Instead, following~\cite{C061},
we take representation~\eqref{crm1} as the
basis and extend it to~$L^0$.

\begin{Definition}\rm
\label{CRM3}
A \textit{coherent utility function on $L^0$} is a map
$u\colon L^0\to[-\infty,\infty]$ defined as
\begin{equation}
\label{crm2}
u(X)=\inf_{\QQ\in\DDD}\EE_\QQ X,\quad X\in L^0,
\end{equation}
where $\DDD$ is a non-empty subset of $\PPP$ and
$\EE_\QQ X$ is understood as $\EE_\QQ X^+-\EE_\QQ X^-$
($X^+=\max\{X,0\}$, $X^-=\max\{-X,0\}$)
with the convention $\infty-\infty=-\infty$.
(Throughout the paper, all the expectations are
understood in this way.)
\end{Definition}

\Remark
This way of defining coherent utility
has parallels to what is done with the classical
expected utility. Namely, the Von Neumann--Morgenstern
representation shows that (appropriately axiomatized)
investor's preferences are described by $\EE U(X)$ with
a bounded function~$U$.
Then one typically takes a concave increasing unbounded
function $U\colon\R\to\R$ and \textit{defines} the preferences
by $\EE U(X)$.\End

\skm
Clearly, different sets $\DDD$ might define the
same coherent utility (for example, $\DDD$ and its
convex hull define the same function~$u$).
However, among all the sets $\DDD$ defining the same~$u$
there exists the largest one.
It is given by
$\{\QQ\in\PPP:\EE_\QQ X\ge u(X)\text{ for any }X\}$.

\begin{Definition}\rm
\label{CRM4}
We will call the largest set, for which~\eqref{crm1}
{\rm(}resp., \eqref{crm2}{\rm)} is true, the
\textit{determining set} of~$u$.
\end{Definition}

\Remarks
(i) Let $\rho'$ be another coherent utility
with the determining set~$\DDD'$. Clearly,
$$
\DDD'\supseteq\DDD\;\Lea\;\rho'\ge\rho.
$$
In other words, the size of $\DDD$ controls the risk
aversion of~$\rho$.

(ii) The determining set is convex.
For coherent utilities on~$L^\infty$, it is also
$L^1$-closed (for the corresponding example,
see~\cite[Subsect.~2.1]{C061}).
In particular, the determining
set of a coherent utility on $L^0$ and the
determining set of its restriction to $L^\infty$ might
be different.

(iii) Let $\DDD$ be an $L^1$-closed convex subset of~$\PPP$.
(Let us note that a particularly important case is
where $\DDD$ is $L^1$-closed, convex, and uniformly
integrable; this condition will be needed in a number
of places below.)
Define a coherent utility~$u$ by~\eqref{crm2}.
Then $\DDD$ is the determining set of~$u$.
Indeed, assume that the determining set $\wt\DDD$ is
larger than~$\DDD$, i.e. there exists
$\QQ_0\in\wt\DDD\setminus\DDD$. Then, by the Hahn--Banach
theorem, we can find $X_0\in L^\infty$ such that
$\EE_{\QQ_0}X_0<\inf_{\QQ\in\DDD}\EE_\QQ X$, which is a
contradiction.
The same argument shows that $\DDD$ is also the determining
set of the restriction of~$u$ to $L^\infty$.\End

\skm
In what follows, we will always consider coherent utility
functions on~$L^0$.

\skb
\textbf{2. Examples.}
Let us now provide several natural examples of coherent
risk measures.

\begin{Example}[Tail V@R]\rm
\label{CRM5}
Let $\la\in(0,1]$ and consider
\begin{equation}
\label{crm3}
\DDD_\la=\Bigl\{\QQ\in\PPP:\frac{d\QQ}{d\PP}
\le\la^{-1}\Bigr\}.
\end{equation}
The corresponding coherent risk measure is called
\textit{Tail V@R} (the terms \textit{Average V@R},
\textit{Conditional V@R}, \textit{Expected Shortfall},
and \textit{Expected Tail Loss} are also used).
Let us denote it by~$\rho_\la$ and the corresponding
coherent utility by~$u_\la$.

Clearly,
$$
\la'\le\la\;\Lea\;\rho_{\la'}\ge\rho_\la,
$$
so that $\la$ serves as the risk aversion parameter.
We have
$$
\rho_\la(X)\xra[\la\da0]{}-\essinf_\omega X(\omega)
$$
(recall that
$\essinf_\omega X(\omega):=\sup\{x:X\ge x\text{ a.s.}\}$).
The right-hand side of this relation is the most severe
risk measure (it is easy to see that any coherent risk
measure~$\rho$ satisfies the inequality
$\rho(X)\le-\essinf_\omega X(\omega)$).
Furthermore, $\rho_1(X)=-\EE X$, which is the most
liberal risk measure (it is seen from Theorem~A.1
that any law invariant risk measure~$\rho$ satisfies
the inequality $\rho(X)\ge-\EE X$).

Let us provide a more explicit representation of~$u_\la$.
Set
$$
Z_*=\la^{-1}I(X<q_\la(X))
+\frac{1-\la^{-1}\PP(X<q_\la(X))}{\PP(X=q_\la(X))}\,
I(X=q_\la(X)).
$$
Throughout the paper, $q_\la$ will denote the
right $\la$-quantile, i.e.
$q_\la(X)=\inf\{x:F(x)\ge\la\}$,
where $F$ is the distribution function of~$X$.
Then $Z_*\ge0$, $\EE Z_*=1$ and,
for any $Z\in\DDD_\la$, we have
\begin{align*}
\EE XZ-\EE XZ_*
&=\EE(X-q_\la(X))Z-\EE(X-q_\la(X))Z_*\\
&=\EE[(X-q_\la(X))(Z-\la^{-1})I(X-q_\la(X)<0)\\
&\hspace*{8mm}+(X-q_\la(X))ZI(X-q_\la(X)>0)]\ge0.
\end{align*}
Hence,
\begin{equation}
\label{crm4}
u_\la(X)=\EE XZ_*
=\la^{-1}\int_{(-\infty,q_\la(X))}x\QQ(dx)
+(1-\la^{-1}\PP(X<q_\la(X))q_\la(X),
\end{equation}
where $\QQ=\Law X$. In particular, it follows that $u_\la$
is law invariant, i.e. it depends only on the distribution
of~$X$.
If $X$ has a continuous distribution, then this formula
is simplified to:
$$
u_\la(X)=\EE(X\cond X\le q_\la(X)).
$$

Using~\eqref{crm4}, one easily gets an equivalent
representation of $u_\la$:
\begin{equation}
\label{crm5}
u_\la(X)=\la^{-1}\int_0^\la q_x(X)dx.
\end{equation}
It is seen from this representation that
$\rho_\la(X)\ge\text{V@R}_\la(X)$.
The advantage of Tail V@R over V@R is that it takes into
account the heaviness of the $\la$-tail (see Figure~1).
Kusuoka~\cite{K01} proved that
on $L^\infty$, $\rho_\la$ is the smallest
law invariant coherent risk measure
that dominates $\text{V@R}_\la$
(the proof can also be found in~\cite[Th.~4.61]{FS04}
or~\cite[Th.~1.48]{S05}).
This suggests an opinion that Tail V@R might be the
most important subclass of coherent risk measures.
However, there exists a risk measure, which is in our
opinion much better than Tail V@R. This is the risk
measure of the next example.\End
\end{Example}

\begin{figure}[!h]
\begin{picture}(150,60)(-35,-27.5)
\put(-10,-0.2){\includegraphics{factor.1}}
\put(-15,0){\vector(1,0){60}}
\put(15,0){\vector(0,1){30}}
\put(55,0){\vector(1,0){50}}
\put(75,0){\vector(0,1){30}}
\put(10.2,0){\line(0,1){2}}
\put(70.2,0){\line(0,1){5}}
\put(8,-4){\small $q_\la$}
\put(68,-4){\small $q_\la$}
\put(-12,-20){\parbox{114mm}{\small\textbf{Figure~1.}
These two distributions have the same
$\la$-quantiles ($\la$ is fixed), so that
$\text{V@R}_\la$ is the same
for them. However, the distribution at the right
is clearly ``better'' than the distribution at the left.}}
\end{picture}
\end{figure}

\begin{Example}[Weighted V@R]\rm
\label{CRM6}
Let $\mu$ be a probability measure on $(0,1]$.
\textit{Weighted V@R with the weighting measure~$\mu$}
(the term \textit{spectral risk measure} is also used)
is the coherent risk measure corresponding to the
coherent utility function
$$
u_\mu(X)=\int_{(0,1]}u_\la(X)\mu(d\la),
$$
where $\int f(x)\mu(dx)$ is understood as
$\int f^+(x)\mu(dx)-\int f^-(x)\mu(dx)$
with the convention $\infty-\infty=-\infty$.
(Throughout the paper, all the integrals are understood
in this way.)
One can check that $u_\mu$ is indeed a coherent utility
(see~\cite[Sect.~3]{C05e} for details).
The measure $\mu$ reflects the risk aversion of~$\rho_\mu$:
the more is the mass of $\mu$ attributed to the left
part of $(0,1]$, the more risk averse is~$\rho_\mu$.

Let us give two arguments in favor of Weighted V@R over
Tail V@R:
\begin{mitemize}
\item (Financial argument) Tail V@R of order $\la$
takes into consideration
only the \mbox{$\la$-tail} of the distribution of $X$; thus,
two distributions with the same $\la$-tail will be
assessed by this measure in the same way, although one
of them might be clearly better than the other (see
Figure~2). On the other hand, if the right endpoint
of the support of~$\mu$ is~$1$, then $\rho_\mu$ depends on the
whole distribution of~$X$.
The paper~\cite{Dowd05} provides some further financial
arguments in favor of Weighted V@R.
\item (Mathematical argument) If the support of~$\mu$
is the whole $[0,1]$, then $\rho_\mu$ possesses
some nice properties that are
not shared by $\rho_\la$. In particular, if $X$ and $Y$
are not comonotone (for the definition, see
Section~\ref{RC}),
then $\rho_\mu(X+Y)<\rho_\mu(X)+\rho_\mu(Y)$
(the proof can be found in~\cite[Sect.~5]{C05e}, where
this was called the \textit{strict diversification
property}). This property is important because it leads
to the uniqueness of a solution of several optimization
problems based on Weighted V@R (see~\cite[Sect.~5]{C05e}).
\end{mitemize}
To put it briefly, Weighted V@R is ``smoother'' than
Tail V@R.

\begin{figure}[!h]
\begin{picture}(150,60)(-26,-27.5)
\put(0,-0.2){\includegraphics{factor.2}}
\put(-5,0){\vector(1,0){50}}
\put(15,0){\vector(0,1){30}}
\put(55,0){\vector(1,0){60}}
\put(75,0){\vector(0,1){30}}
\put(10.2,0){\line(0,1){4}}
\put(70.2,0){\line(0,1){4}}
\put(8,-4){\small $q_\la$}
\put(68,-4){\small $q_\la$}
\put(-2,-20){\parbox{114mm}{\small\textbf{Figure~2.}
These two distributions have the same
$\la$-tails ($\la$ is fixed), so that
$\text{TV@R}_\la$ is the same for
them. However, the distribution at the right is clearly
``better'' than the distribution at the left.}}
\end{picture}
\end{figure}

Weighted V@R is law invariant because Tail V@R possesses
this property. Kusuoka~\cite{K01} proved that
on $L^\infty$ the class of Weighted V@Rs is exactly the
class of law invariant coherent risk measures satisfying
the additional property of \textit{comonotonicity}, which
means that $\rho(X+Y)=\rho(X)+\rho(Y)$ for any
comonotone (the definition is recalled in Section~\ref{RC})
random variables $X$ and~$Y$
(the proof can also be found in~\cite[Th.~4.87]{FS04}
or~\cite[Th.~1.58]{S05}).

Let us now provide several equivalent representations
of Weighted V@R. It follows from~\eqref{crm5} that
\begin{equation}
\label{crm6}
u_\mu(X)
=\int_{(0,1]}\la^{-1}\int_0^\la q_x(X)dx\mu(d\la)\\
=\int_0^1 q_x(X)\psi_\mu(x)dx,
\end{equation}
where
\begin{equation}
\label{crm7}
\psi_\mu(x)=\int_{[x,1]}\la^{-1}\mu(d\la),
\quad x\in[0,1].
\end{equation}
The last formula establishes a one-to-one correspondence
between the left-continuous decreasing functions
$\psi\colon[0,1]\to[0,1]$ with $\int_{(0,1]}\psi(x)dx=1$
and the probability measures on $(0,1]$.

In particular, let $\Omega=\{1,\dots,T\}$ and
$X(t)=x_t$. Let $x_{(1)},\dots,x_{(T)}$ be the values
$x_1,\dots,x_T$ in the increasing order.
Define $n(t)$ through the equality $x_{(t)}=x_{n(t)}$.
Then
\begin{equation}
\label{crm8}
u_\mu(X)=\sum_{t=1}^T x_{n(t)}
\int_{z_{t-1}}^{z_t}\psi_\mu(x)dx,
\end{equation}
where $z_t=\sum_{i=1}^t\PP\{n(i)\}$.
This formula yields a simple procedure of the
empirical estimation of~$u_\mu(X)$.

In order to provide another representation of $u_\mu$,
consider the function
\begin{equation}
\label{crm9}
\Psi_\mu(x)
=\int_0^x\psi_\mu(y)dy
=\int_0^x\int_{(y,1]}\la^{-1}\mu(d\la)dy,
\quad x\in(0,1].
\end{equation}
It is easy to see that $\Psi_\mu\colon[0,1]\to[0,1]$ is
increasing, concave, continuous, $\Psi_\mu(0)=0$, and
$\Psi_\mu(1)=1$.
In fact, \eqref{crm9} establishes a one-to-one correspondence
between the functions with these properties and the
probability measures~$\mu$ on $(0,1]$
(for details, see~\cite[Lem.~4.63]{FS04} or~\cite[Lem.~1.50]{S05}).
The inverse map $\Psi\mapsto\mu$ is given by
$\mu=-\la\Psi''$, where $\Psi''$ is the second
derivative of $\Psi$ taken in the sense of distributions,
i.e. it is the measure on $(0,1]$ defined by
$\Psi''((a,b]):=\Psi'_+(b)-\Psi'_+(a)$, where $\Psi'_+$
is the right-hand derivative.
Let $F$ be the distribution function of~$X$.
As the function $x\mapsto q_x(X)$ is constant on the
intervals of the form $[F(y-),F(y))$, we can derive
from~\eqref{crm5} the following representation:
\begin{equation}
\label{crm10}
u_\mu(X)
=\int_0^1 q_x(X)d\Psi_\mu(x)
=\int_\R q_{F(x)}(X)d\Psi_\mu(F(x))
=\int_\R xd\Psi_\mu(F(x))
=\EE Y,
\end{equation}
where $Y$ is a random variable with the distribution
function $\Psi_\mu\circ F$.
This representation provides a convenient tool
for designing particular risk measures.
Let us remark that the functionals of the form~\eqref{crm10}
were considered by actuaries under the name
\textit{distorted measures}
already in the early 90s, i.e. before the papers of
Artzner, Delbaen, Eber, and Heath; see, for
example,~\cite{D90}, \cite{W96}
(see also the paper~\cite{WYP97}, which appeared at the
same time as~\cite{ADEH97}).

If $X$ has a continuous distribution, we get from~\eqref{crm10}
one more representation:
\begin{equation}
\label{crm11}
u_\mu(X)
=\int_\R x\psi_\mu(F(x))dF(x)
=\EE X\psi_\mu(X)
=\EE_{\QQ_\mu(X)}X,
\end{equation}
where $\QQ_\mu(X)=\psi_\mu(X)\PP$.
Note that
$$
\EE\psi_\mu(F(X))
=\int_\R\psi_\mu(F(x))dF(x)
=\int_0^1\psi_\mu(x)dx
=\int_{(0,1]}\int_0^\la\la^{-1}dx\mu(d\la)
=1,
$$
so that $\QQ_\mu(X)$ is a probability measure.
As $\psi_\mu\circ F$ is decreasing, this measure attributes
more mass to the outcomes corresponding to low values
of~$X$. Thus, $\QQ_\mu(X)$ reflects the risk aversion of
an agent who possesses the position that yields the
P\&L~$X$.

The determining set $\DDD_\mu$ of $u_\mu$ admits the
following representation:
\begin{equation}
\label{crm12}
\DDD_\mu=\{Z\in L^0:Z\ge0,\;\EE Z=1,\;\text{\rm and }
\EE(Z-x)^+\le\Phi_\mu(x)\;\forall x\in\R_+\},
\end{equation}
where
\begin{equation}
\label{crm13}
\Phi_\mu(x)=\sup_{y\in[0,1]}(\Psi_\mu(y)-xy),\quad x\in\R_+
\end{equation}
(see~\cite[Th.~4.6]{C05e}).
Note that $\Phi_\mu\colon\R_+\to\R_+$ is continuous, decreasing,
convex, $\Phi_\mu(0)=1$, $\Phi_\mu(x)>0$ for
$x<\int_{(0,1]}\la^{-1}d\mu$, and $\Phi_\mu(x)=0$ for
$x\ge\int_{(0,1]}\la^{-1}d\mu$.

\begin{figure}[!h]
\begin{picture}(150,65)(-18,-17.5)
\put(-0.4,-0.3){\includegraphics{factor.3}}
\put(0,0){\vector(1,0){25}}
\put(0,0){\vector(0,1){45}}
\put(-0.5,40){\line(1,0){1}}
\put(20,-0.5){\line(0,1){1}}
\put(-15,39){\scalebox{0.7}{$\int_{(0,1]}\la^{-1}d\mu$}}
\put(19,-4){\small$1$}
\put(-1,-4){\small$0$}
\put(24,-3){\small $x$}
\put(4,17){\small$\psi_\mu$}
\put(35,0){\vector(1,0){25}}
\put(35,0){\vector(0,1){25}}
\multiput(55,0.5)(0,2){10}{\line(0,1){1}}
\multiput(35.5,20)(2,0){10}{\line(1,0){1}}
\put(32,19){\small$1$}
\put(34,-4){\small$0$}
\put(54,-4){\small$1$}
\put(59,-3){\small $x$}
\put(37.5,3.5){\scalebox{0.7}{$\int_{(0,1]}\la^{-1}d\mu$}}
\put(36,15){\small$\Psi_\mu$}
\put(70,0){\vector(1,0){55}}
\put(70,0){\vector(0,1){25}}
\put(69.5,20){\line(1,0){1}}
\put(110,-0.5){\line(0,1){1}}
\put(67,19){\small$1$}
\put(69,-4){\small$0$}
\put(124,-3){\small $x$}
\put(103,-4){\scalebox{0.75}{$\int_{(0,1]}\la^{-1}d\mu$}}
\put(84,10){\small$\Phi_\mu$}
\put(28,-15){\small\textbf{Figure~3.} The form of
$\psi_\mu$, $\Psi_\mu$, and $\Phi_\mu$}
\end{picture}
\end{figure}

One more representation of $\DDD_\mu$ is:
\begin{align*}
\DDD_\mu
&=\{\QQ\in\PPP:\QQ(A)\le\Psi_\mu(\PP(A))\;
\text{\rm for any }A\in\F\}\\
&=\Bigl\{Z\in L^0:Z\ge0,\;\EE_\PP Z=1,\text{ and }
\int_{1-x}^1\!\!q_s(Z)ds\le\Psi_\mu(x)\;\forall
x\in[0,1]\Bigr\}.
\end{align*}
It was obtained in~\cite{CD03}
(the proof can also be found in~\cite[Th.~4.73]{FS04}
or~\cite[Th.~1.53]{S05}).
It is seen from this representation that
\begin{equation}
\label{crm14}
\Psi_{\mu'}\ge\Psi_\mu\;\Lea\;\rho_{\mu'}\ge\rho_\mu.
\end{equation}

Moreover, it is easy to see that
\begin{equation}
\label{crm15}
\mu'\preccurlyeq\mu\;\Longrightarrow\;
\rho_{\mu'}\ge\rho_\mu,
\end{equation}
where the notation $\mu'\preccurlyeq\mu$ means that
$\mu$ stochastically dominates $\mu'$, i.e. their
distribution functions satisfy $F_{\mu'}\ge F_\mu$.
In order to prove~\eqref{crm15}, it is sufficient to
notice that $u_\la$ is increasing in~$\la$ and
$u_\mu(X)=\wt\EE u_\xi(X)$, where
$\xi$ is a random variable on some space
$(\wt\Omega,\wt\F,\wt\PP)$ with $\Law\xi=\mu$.
One should also use the following well-known fact
(see~\cite[\S~1.A]{ShSh94}): $\mu'\preccurlyeq\mu$
if and only on some space there exist $\xi'\le\xi$
with $\Law\xi=\mu$, $\Law\xi'=\mu'$.

For more information on Weighted V@R, see~\cite{A02},
\cite{A04}, \cite{C05e}, \cite{Dowd05}.\End
\end{Example}

Weighted V@R is a particular case of the more general
class of risk measures that is described in the next
example.

\begin{Example}[Law invariant risk measures]\rm
\label{CRM7}
A risk measure $\rho$ is \textit{law invariant} if
$\rho(X)=\rho(Y)$ whenever $\Law X=\Law Y$.
Kusuoka~\cite{K01} proved that
a coherent risk measure $\rho$ on $L^\infty$
is law invariant if and only if it has the form
$$
\rho(X)=\sup_{\mu\in\mathfrak M}\rho_\mu(X)
$$
with some set $\mathfrak M$ of probability measures
on $(0,1]$ (the proof can also be found
in~\cite[Cor.~4.58]{FS04} or~\cite[Cor.~1.45]{S05}).
Theorem~A.1 extends this result to risk measures
on~$L^0$.\End
\end{Example}

An interesting two-parameter family of coherent risk
measures is provided by the example below.

\begin{Example}[Moment-based risk measures]\rm
\label{CRM8}
Let $p\in[1,\infty]$ and $\al\in[0,1]$.
Consider the set
$$
\DDD=\{1+\al(Z-\EE Z):Z\ge0,\,\|Z\|_q\le1\},
$$
where $q=p/(p-1)$ and $\|Z\|_q=(\EE Z^q)^{1/q}$.
Then the corresponding coherent risk measure has the form
$$
\rho(X)=-\EE X+\al\|(X-\EE X)^-\|_p
$$
(see~\cite[Sect.~4]{D05}).
In particular, if $p=2$, $\al=1$, and $\EE X=0$, then
$\rho(X)$ is the semivariance of~$X$. Semivariance was
proposed by Markowitz~\cite{M59} as a substitute for variance.
Its advantage is that it measures really risk
(i.e. the downfall); its disadvantage is that it is less
convenient analytically than variance.\End
\end{Example}

The moment-based risk measures are law invariant, but they do not
belong to the class of Weighted V@Rs, which is a very
convenient class.
Weighted V@Rs are parametrized by the probability
measures~$\mu$ on $(0,1]$, which is a huge class.
For the practical purposes, one needs to select a
convenient finite-parameter subclass of Weighted V@Rs.
The most natural parametric family of probability
measures on $[0,1]$ is the family of Beta distributions.
It turns out that the family of corresponding Weighted
V@Rs admits a very natural interpretation and a very
simple estimation procedure.
We call these measures Beta V@Rs in accordance with
Beta distributions.

\begin{Example}[Beta V@R]\rm
\label{CRM9}
Let $\al\in(-1,\infty)$, $\be\in(-1,\al)$.
\textit{Beta V@R with parameters $\al,\be$} is the
Weighted V@R with the weighting measure
$$
\mu_{\al,\be}(dx)={\rm B}(\be+1,\al-\be)^{-1}x^\be(1-x)^{\al-\be-1}dx,
\quad x\in[0,1].
$$
This risk measure will be denoted as $\rho_{\al,\be}$
and the corresponding coherent utility will be denoted
as $u_{\al,\be}$.

It follows from~\eqref{crm15} that
\begin{align}
\label{crm16}
&\al'\ge\al\;\Lea\;
\mu_{\al',\be}\preccurlyeq\mu_{\al,\be}\;\Lea\;
\rho_{\al',\be}\ge\rho_{\al,\be},\\
\label{crm17}
&\be'\ge\be\;\Lea\;
\mu_{\al,\be'}\succcurlyeq\mu_{\al,\be}\;\Lea\;
\rho_{\al,\be'}\le\rho_{\al,\be}.
\end{align}
Furthermore,
\begin{align*}
\mu_{\al,\be}\xra[\al\to\infty]{}\de_0
&\;\Longrightarrow\;
\rho_{\al,\be}(X)\xra[\al\to\infty]{}
-\essinf_\omega X(\omega),\\
\mu_{\al,\be}\xra[\be\ua\al]{}\de_1
&\;\Longrightarrow\;
\rho_{\al,\be}(X)\xra[\be\ua\al]{}-\EE X,
\end{align*}
where $\de_a$ denotes the delta-mass concentrated at~$a$.
In particular, these relations show that we can redefine
$\rho_{\al,\al}(X)$ as $-\EE X$.

Suppose now that $\al,\be\in\N$.
Let $X_1,\dots,X_\al$ be independent copies of~$X$,
$X_{(1)},\dots,X_{(\al)}$ be the corresponding order
statistics, and $\xi$ be an independent uniformly
distributed on $\{1,\dots,\be\}$ random variable.
Let us prove that
\begin{equation}
\label{crm18}
\rho_{\al,\be}(X)
=-\EE X_{(\xi)}
=-\EE\Bigl[\frac{1}{\be}\sum_{i=1}^\be X_{(i)}\Bigr].
\end{equation}
The second equality is obvious, so we should prove only
the first one.
The random variables $X_i$ can be realized as
$X_i=F^{-1}(U_i)$, where $F$ is the distribution
function of~$X$ and $U_1,\dots,U_\al$ are independent
uniformly distributed on $[0,1]$ random variables.
Let $U_{(1)},\dots,U_{(\al)}$ denote their order
statistics. We have (see~\cite[Ch.~1, \S~7]{F71})
$$
\frac{d}{dx}\,\PP(U_{(i)}\le x)
=\al C_{\al-1}^{i-1}x^{i-1}(1-x)^{\al-i},
\quad x\in(0,1),
$$
so that
$$
\frac{d}{dx}\,\PP(U_{(\xi)}\le x)
=\frac{1}{\be}\sum_{i=1}^\be\frac{d}{dx}\,
\PP(U_{(i)}\le x)\\
=\frac{\al}{\be}\sum_{i=1}^{\be}
C_{\al-1}^{i-1}x^{i-1}(1-x)^{\al-i},\quad x\in(0,1)
$$
and
$$
\frac{d^2}{dx^2}\,\PP(U_{(\xi)}\le x)
=-\frac{\al!}{\be!(\al-\be-1)!}x^{\be-1}(1-x)^{\al-\be-1},
\quad x\in(0,1).
$$
For the function
$\Psi_{\al,\be}:=\Psi_{\mu_{\al,\be}}$, we have
$$
\Psi''_{\al,\be}(x)
=-x^{-1}{\rm B}(\be+1,\al-\be)^{-1}x^\be(1-x)^{\al-\be-1}dx
=\frac{d^2}{dx^2}\,\PP(U_{(\xi)}\le x),
\quad x\in(0,1).
$$
Moreover, the functions $\Psi_{\al,\be}$ and
$\PP(U_{(\xi)}\le\cdot)$ coincide at~0 and at~1.
Consequently, these functions are equal, and therefore,
$$
\PP(X_{(\xi)}\le x)
=\PP(F^{-1}(U_{(\xi)})\le x)
=\PP(U_{(\xi)}\le F(x))
=\Psi_{\al,\be}(F(x)),\quad x\in\R.
$$
Recalling~\eqref{crm10}, we obtain~\eqref{crm18}.

It is seen from the calculations given above that for
Beta V@R the function
$\psi_{\al,\be}:=\psi_{\mu_{\al,\be}}$ has the form
\begin{equation}
\label{crm19}
\psi_{\al,\be}(x)
=\Psi'_{\al,\be}(x)
=\frac{d}{dx}\,\PP(U_{(\xi)}\le x)
=\frac{\al}{\be}\sum_{i=1}^\be C_{\al-1}^{i-1}
x^{i-1}(1-x)^{\al-i},\quad x\in(0,1).
\end{equation}

Let us remark that, according to~\eqref{crm5}, Tail V@R
admits the following representation:
$\rho_\la(X)=\wt\EE q_\xi(X)$, where $\xi$ is a random
variable on some space $(\wt\Omega,\wt\F,\wt\PP)$ with
the uniform distribution on $[0,\la]$.
Formula~\eqref{crm18} can be rewritten as follows:
$\rho_{\al,\be}(X)=-\wt\EE q_\xi(\wt X)$, where
$\wt X$ is a random variable distributed according
to the empirical distribution constructed
by $X_1,\dots,X_\al$,
$\xi$ is an independent random variable with the
uniform distribution on $[0,\be/\al]$,
and $\wt\EE$ means the averaging
over empirical distributions.\End
\end{Example}

An important one-parameter family of coherent risk
measures is obtained from Beta V@R by fixing the
value $\be=1$.

\begin{Example}[Alpha V@R]\rm
\label{CRM10}
Let $\al\in(1,\infty)$.
\textit{Alpha V@R of order $\al$} is Beta V@R of order
$(\al,1)$.
It is seen from~\eqref{crm16} that $\al$ measures the
risk aversion of~$\rho_\al$. Furthermore, it follows
from~\eqref{crm18} that, for $\al\in\N$,
\begin{equation}
\label{crm20}
\rho_\al(X)=-\EE\min_{i=1,\dots,\al}X_i,
\end{equation}
where $X_1,\dots,X_\al$ are independent copies of~$X$.\End
\end{Example}

The classes of risk measures described by
Examples~\ref{CRM5}--\ref{CRM10} are related by the
following diagram:
\begin{center}
\begin{tabular}{ccccccc}
\raisebox{-2mm}{\vbox{
\hbox{Tail}
\vspace{2mm}
\hbox{\hspace*{-1mm}V@R}}}&
$\subset$&
\raisebox{-2mm}{\vbox{
\hbox{Weighted}
\vspace{2mm}
\hbox{\hspace*{3.5mm}V@R}}}&
$\subset$&
\raisebox{-2mm}{\vbox{
\hbox{Law-invariant}
\vspace{2mm}
\hbox{risk measures}}}&
$\subset$&
\raisebox{-2mm}{\vbox{
\hbox{\hspace*{4mm}Coherent}
\vspace{2mm}
\hbox{risk measures}}}\\[4mm]
&&\rotatebox{90}{$\subset$}&&
\rotatebox{90}{$\subset$}\\[1mm]
&&\raisebox{-2mm}{\vbox{
\hbox{Beta}
\vspace{2mm}
\hbox{\hspace*{-0.2mm}V@R}}}&&
\raisebox{-2mm}{\vbox{
\hbox{Moment-based}
\vspace{2mm}
\hbox{\hspace*{1mm}risk measures}}}&&\\[4mm]
&&\rotatebox{90}{$\subset$}&&&&\\[1mm]
&&\raisebox{-5mm}{\vbox{
\hbox{Alpha}
\vspace{2mm}
\hbox{\hspace*{1mm}V@R}}}&&&&\\
\end{tabular}
\end{center}
In our opinion, the best classes are: Weighted V@R,
Beta V@R, and Alpha V@R. All the empirical estimation
procedures considered in the paper will be provided
for these three classes.

\skm
\textbf{3. Further examples.}
Coherent risk measures are primarily intended to assess
the risk of non-Gaussian P\&Ls. However, as an example,
it is interesting to look at their values in
the Gaussian case.

\begin{Example}\rm
\label{CRM11}
\textbf{(i)} Let $u$ be a law invariant coherent utility
that is finite on Gaussian random variables.
It is easy to see that then there exists $\ga\in\R_+$
such that, for any Gaussian random variable~$X$ with
mean~$m$ and variance~$\si^2$, we have
$u(X)=m-\ga\si$.
In particular, if $X$ is a $d$-dimensional Gaussian
random vector with mean~0 and covariance matrix~$C$,
then $\rho(\lb h,X\rb)=\ga\lb h,Ch\rb$, $h\in\R^d$.

\textbf{(ii)} For Tail V@R, we get an explicit form
of the constant~$\ga$:
$$
\ga(\la)
=-\frac{1}{\la\sqrt{2\pi}}\int_{-\infty}^{q_\la}x e^{-x^2/2}dx
=\frac{1}{\la\sqrt{2\pi}}\int_{q_\la^2/2}^\infty e^{-y}dy
=\frac{1}{\la\sqrt{2\pi}}e^{-q_\la^2/2},
$$
where $q_\la$ is the $\la$-quantile of the standard normal
distribution (in order to check the second equality,
one should consider separately the cases $\la\le1/2$
and $\la\ge1/2$).

\begin{picture}(150,67.5)(-34,-15)
\put(-7.5,-4.4){\includegraphics{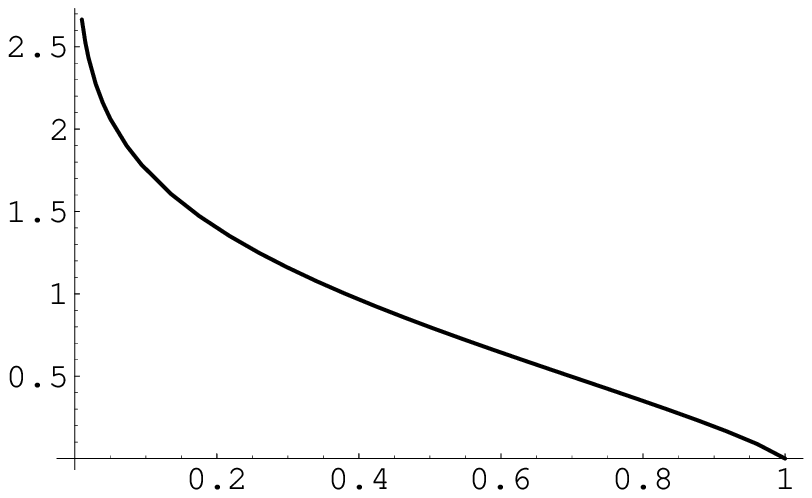}}
\put(0,0){\vector(1,0){80}}
\put(0,0){\vector(0,1){50}}
\put(79,-3){\small$\la$}
\put(25,20){\small$\ga(\la)$}
\put(15,-12){\small\textbf{Figure~4.} The form of $\ga(\la)$}
\end{picture}

\textbf{(iii)} For Beta V@R, we have from~(ii):
$$
\ga(\al,\be)=\int_0^1\frac{1}{\sqrt{2\pi}\,
{\rm B}(\be+1,\al-\be)}\,
e^{-q_x^2/2}x^{\be-1/2}(1-x)^{\al-\be-1}dx.
$$

\begin{picture}(150,67.5)(-34,-15)
\put(-7.7,-104.5){\scalebox{0.8}{\includegraphics{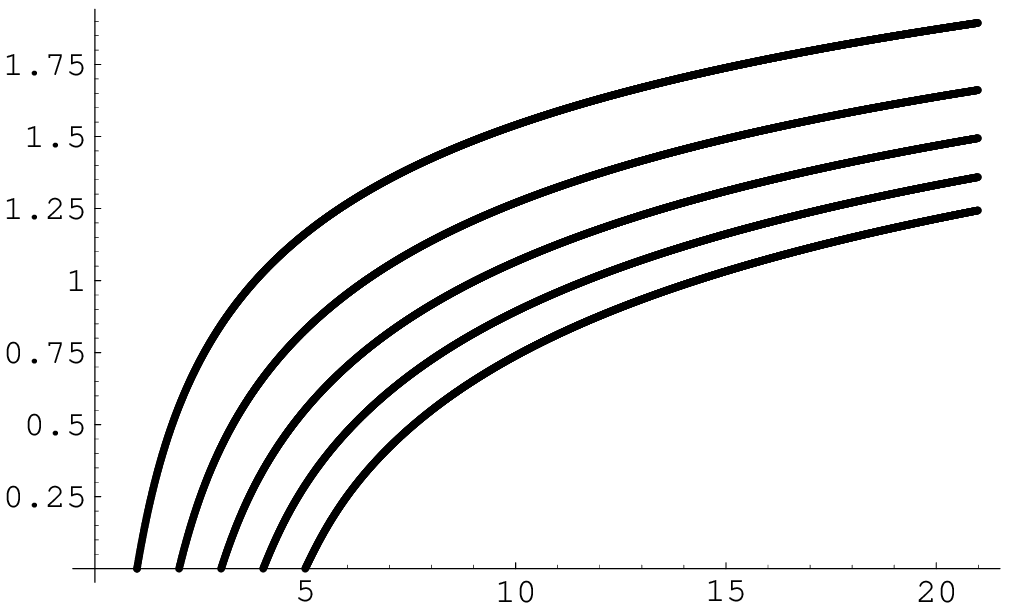}}}
\put(0,0){\vector(1,0){80}}
\put(0,0){\vector(0,1){50}}
\put(79,-3){\small $\al$}
\put(72,27){\scalebox{0.8}{$\ga(\al,5)$}}
\put(72,30.5){\scalebox{0.8}{$\ga(\al,4)$}}
\put(72,34){\scalebox{0.8}{$\ga(\al,3)$}}
\put(72,38){\scalebox{0.8}{$\ga(\al,2)$}}
\put(72,43.5){\scalebox{0.8}{$\ga(\al,1)$}}
\put(14,-12){\small\textbf{Figure~5.} The form of
$\ga(\al,\be)$}
\end{picture}
\end{Example}

Let us also give a nice credit risk example,
which illustrates the effect of coherent risk
diversification achieved in large portfolios.
The example is borrowed from~\cite{D05}.

\begin{Example}\rm
\label{CRM12}
Let $\mu$ be a measure on $(0,1]$ such that
$\int_{(0,1]}\la^{-1}\mu(d\la)<\infty$ (for example,
the weighting measures of Tail V@R, Beta V@R with
$\be>0$, and Alpha V@R satisfy this condition).
Let $X_1,X_2,\dots$ be independent identically distributed
integrable random variables.
Set $S_n=X_1+\dots+X_n$. By the law of large numbers,
$S_n/n\xra{L^1}\EE X$.
Consequently, $u_\la(S_n/n)\to\EE X$ for any $\la\in(0,1]$.
Using the estimate $|u_\la(\xi)|\le\la^{-1}\EE|\xi|$
and the Lebesgue dominated convergence theorem, we get
$$
u_\mu\Bigl(\frac{S_n}{n}\Bigr)
=\int_{(0,1]}u_\la\Bigl(\frac{S_n}{n}\Bigr)\mu(d\la)
\xra[n\to\infty]{}\int_{(0,1]}u_\la(\EE X)\mu(d\la)
=\EE X.
$$
This result admits the following interpretation.
If a firm gives many loans of size~$L$ to independent
identical customers (so that the $n$-th
customer returns back the random amount~$X_n$),
then the firm is on the safe side provided that
$\EE X\ge L$.\End
\end{Example}

\textbf{4. Empirical estimation.}
Here we will describe procedures for the empirical
estimation of Alpha V@R, Beta V@R, and Weighted V@R,
which serve as coherent counterparts of the historic
V@R estimation (see~\cite[Sect.~6]{M02}).
Let $X$ be the increment of the value of some portfolio
over the unit time period~$\De$.

In order to estimate Alpha V@R with $\al\in\N$,
one should first choose the number of trials $K\in\N$
and generate independent draws
$(x_{kl};k=1,\dots,K,\,l=1,\dots,\al)$ of~$X$.
This can be done by one of the following techniques:
\begin{mitemize}
\item Each $x_{kl}$ is drawn uniformly from the recent
$T$ realizations $x_1,\dots,x_T$ of~$X$.
For example, if $\De$ is one day (which is a typical choice),
these are $T$ recent daily increments of the value of
the portfolio under consideration.
\item Each $x_{kl}$ is drawn from recent $T$ realizations
$x_1,\dots,x_T$ of~$X$, according to a probability
measure~$\nu$ on $\{1,\dots,T\}$. A natural example is:
$T=\infty$, so that $x_t$ is the increment of the
portfolio's value over the interval $[-t\De,-(t-1)\De]$,
where 0 is the current time instant;
$\nu$ is the geometric distribution with a parameter~$\la$.
This method enables one to put more mass on recent
realizations of~$X$.
It is at the basis of the \textit{weighted historical
simulation} (see~\cite[Sect.~5.3]{C03}).
Typically, $\la$ is chosen between 0.95 and 0.99.
\item Each $x_{kl}$ might be generated using the
bootstrap method, i.e. we split the time axis into
small intervals of length $n^{-1}\De$ and create each
$x_{kl}$ as a sum of the increments of the portfolio's
value over $n$ randomly chosen small intervals.
The bootstrap method might be combined with the
weighting described above, i.e. we take recent small
intervals with a higher probability than old ones.
\end{mitemize}
Having generated $x_{kl}$, one should calculate the
array
$$
l_k=\argmin_{l=1,\dots,\al}x_{kl},\quad k=1,\dots,K.
$$
According to~\eqref{crm20}, an empirical estimate of
$\rho_\al(X)$ is provided by
$$
\rho_\text{e}(\wt X)=-\frac{1}{K}\sum_{k=1}^K x_{kl_k}.
$$

In order to estimate Beta V@R with $\al,\be\in\N$,
one should generate $x_{kl}$ similarly.
Let $l_{k1},\dots,l_{k\be}$ be the numbers
$l\in\{1,\dots,\al\}$ such that the corresponding
$x_{kl}$ stand at the first $\be$ places
(in the increasing order) among
$x_{k1},\dots,x_{k\al}$.
According to~\eqref{crm18}, an empirical estimate of
$\rho_{\al,\be}(X)$ is provided by
$$
\rho_{\text{e}}(X)
=-\frac{1}{K\be}\sum_{k=1}^K\sum_{i=1}^\be x_{kl_{ki}}.
$$

In order to estimate Weighted V@R, one should fix a
data set $x_1,\dots,x_T$ and a measure~$\nu$ on
$\{1,\dots,T\}$.
This can be done by one of the following techniques:
\begin{mitemize}
\item The values $x_1,\dots,x_T$ are recent $T$ realizations
of~$X$.
\item The values $x_1,\dots,x_T$ are obtained through
the bootstrap technique, i.e. each $x_t$ is a sum of
the increments of the portfolio's value over randomly
chosen small intervals; $\nu$ is uniform.
\end{mitemize}
Let $x_{(1)},\dots,x_{(T)}$ be the values $x_1,\dots,x_T$
in the increasing order.
Define $n(t)$ through the equality $x_{(t)}=x_{n(t)}$.
According to~\eqref{crm8}, an empirical estimate of
$\rho_\mu(X)$ is provided by
$$
\rho_\text{e}(X)=
-\sum_{t=1}^T x_{n(t)}\int_{z_{t-1}}^{z_t}\psi_\mu(x)dx,
$$
where $z_t=\sum_{i=1}^T\nu\{n(i)\}$ and $\psi_\mu$
is given by~\eqref{crm7}.

An advantage of Weighted V@R over Alpha V@R and Beta V@R
is that it is a wider class. However, Beta V@R is already
rather a flexible family.
A big advantage of Alpha V@R and Beta V@R is that their
empirical estimation procedure does not require the
data ordering (the number of operations required to
order the data set $x_1,\dots,x_T$ is $T\ln T$;
this is a particularly unpleasant number for
$T=\infty$, which is one of typical possible
choices in estimating Alpha V@R and Beta V@R).
Note that, for estimating V@R from a data set of
size~$T$, one needs to order the time series, and
the number of operations required grows quadratically
in~$T$. Thus, Alpha V@R and Beta V@R are not only much
wiser than V@R; they are also estimated faster!

%=======================================================
\section{Factor Risk}
\label{FR}

\textbf{1. $L^1$-spaces.}
Let $(\Omega,\F,\PP)$ be a probability space
and $u$ be a coherent utility with the
determining set~$\DDD$.

For the theorems below, we need to define the $L^1$-spaces
associated with a coherent risk measure (they were
introduced in~\cite{C061}).
The \textit{weak} and \textit{strong $L^1$-spaces} are
\begin{align*}
L_w^1(\DDD)&=\{X\in L^0:u(X)>-\infty,\,u(-X)>-\infty\},\\
L_s^1(\DDD)&=\Bigl\{X\in L^0:\lim_{n\to\infty}
\sup_{\QQ\in\DDD}\EE_\QQ|X|I(|X|>n)=0\Bigr\}.
\end{align*}
Clearly, $L_s^1(\DDD)\subseteq L_w^1(\DDD)$.
In general, this inclusion might be strict.
Indeed, let $X_0$ be a positive unbounded random variable
with $\PP(X_0=0)>0$ and let $\DDD=\{\QQ\in\PPP:\EE_\QQ X_0=1\}$.
Then $X_0\in L_w^1(\DDD)$, but $X_0\notin L_s^1(\DDD)$.
However, as shown by the examples below, in most
natural situations these two spaces coincide and have
a very simple form.

\begin{Example}\rm
\label{FR1}
\textbf{(i)} If $\DDD=\{\QQ\}$ is a singleton, then
$L_w^1(\DDD)=L_s^1(\DDD)=L^1(\QQ)$, which motivates
the notation.

\textbf{(ii)} For Weighted V@R, we have
$L_w^1(\DDD_\mu)=L_s^1(\DDD_\mu)$
(see~\cite[Subsect.~2.2]{C061}),
so we can denote this space simply by $L^1(\DDD_\mu)$.
It is clear from~\eqref{crm12} that $\PP\in\DDD_\mu$,
so that $L^1\subseteq L^1(\DDD_\mu)$
(throughout the paper, $L^1$ stands for $L^1(\PP)$).
In general, this inclusion can be strict.
However, if $\int_{(0,1]}\la^{-1}\mu(d\la)<\infty$
(for example, the weighting measures of Tail V@R, Beta
V@R with $\be>0$, and Alpha V@R satisfy this condition),
then it is seen from~\eqref{crm12} that all the
densities from~$\DDD_\mu$ are bounded by
$\int_{(0,1]}\la^{-1}\mu(d\la)$, so that
$L^1(\DDD_\mu)=L^1$.

\textbf{(iii)} If $\rho$ is the moment-based risk
measure of Example~\ref{CRM8} with $\al>0$, then clearly
$L_w^1(\DDD)\subseteq L^p$ and
$L^p\subseteq L_s^1(\DDD)$, so that
$L_w^1(\DDD)=L_s^1(\DDD)=L^p$.\End
\end{Example}

\textbf{2. Factor risk.}
Let $Y$ be a random variable (resp., random vector)
meaning the increment over the unit time period
of some market factor (resp., several market factors).

\begin{Definition}\rm
\label{FR2}
The \textit{factor utility} is
$$
u^f(X;Y)=\inf_{\QQ\in\EE(\DDD\cond Y)}\EE_\QQ X,
$$
where $\EE(\DDD\cond Y):=\{\EE(Z\cond Y):Z\in\DDD\}$.

The \textit{factor risk} is $\rho^f(X;Y):=-u^f(X;Y)$.
\end{Definition}

\Remarks
(i) As $\DDD$ is convex, $\EE(\DDD\cond Y)$ is also
convex.

(ii) If $\DDD$ is convex and $L^1$-closed, then
$\EE(\DDD\cond Y)$ is not necessarily $L^1$-closed.
As an example, consider $\Omega=[0,1]^2$ endowed with
the Lebesgue measure and let
$$
\DDD=\Bigl\{\sum_{n=1}^\infty a_nZ_n:a_n\in\R_+,\,
\sum_{n=1}^\infty a_n=1\Bigr\},
$$
where
$$
Z_n(x_1,x_2)=\begin{cases}
1/n&\text{if}\;\;x_1\le1/2,\\
2n-1&\text{if}\;\;x_1>1/2\text{ and }x_2\le1/n,\\
0&\text{otherwise}.
\end{cases}
$$
Let $Y(x^1,x^2)=x^1$. Then
$$
\EE(Z_n\cond Y)\xra[n\to\infty]{L^1}2I(x_1>1/2)
\notin\EE(\DDD\cond Y).
$$

(iii) If $\DDD$ is convex, $L^1$-closed, and uniformly
integrable, then $\EE(\DDD\cond Y)$ is also convex,
$L^1$-closed, and uniformly integrable.
Indeed, let $Z_n\in\DDD$ be such that
$\EE(Z_n\cond Y)\xra{L^1}Z$.
By Komlos' principle of subsequences (see~\cite{K67}),
we can select a sequence $\wt Z_n\in\conv(Z_n,Z_{n+1},\dots)$
such that $\wt Z_n\xra{\text{a.s.}}\wt Z$.
As $\DDD$ is convex and uniformly integrable,
the convergence holds in $L^1$ and $\wt Z\in\DDD$.
Then $\EE(\wt Z\cond Y)=Z$, so that $\EE(\DDD\cond Y)$
is $L^1$-closed.
Its uniform integrability is a well-known fact.\End

\begin{Theorem}
\label{FR3}
For
$X\in L_s^1(\DDD)\cap L_s^1(\EE(\DDD\cond Y))\cap L^1$,
we have
$$
u^f(X;Y)=u(\EE(X\cond Y)).
$$
\end{Theorem}

This theorem follows from

\begin{Lemma}
\label{FR4}
For $X\in L_s^1(\DDD)\cap L_s^1(\EE(\DDD\cond Y))\cap L^1$
and $Z\in\DDD$, we have
$\EE X\EE(Z\cond Y)=\EE\EE(X\cond Y)Z$.
\end{Lemma}

\Proof
We can write $X=\xi_n+\eta_n$, where
$\xi_n=XI(|X|\le n)$, $\eta_n=XI(|X|>n)$, $n\in\N$.
By the definition of $L^1$-spaces,
$\EE|\eta_n|Z\to0$ and $\EE|\eta_n|\EE(Z\cond Y)\to0$.
The equality
$$
\EE\xi_n\EE(Z\cond Y)
=\EE\EE(\xi_n\cond Y)\EE(Z\cond Y)
=\EE\EE(\xi_n\cond Y)Z
$$
and the estimates
$$
\EE|\EE(\eta_n\cond Y)|Z
\le\EE\EE(|\eta_n|\cond Y)Z
=\EE\EE(|\eta_n|\cond Y)\EE(Z\cond Y)
=\EE|\eta_n|\EE(Z\cond Y)
$$
yield the desired statement.\End

\skm
Recall that a probability space $(\Omega,\F,\PP)$ is
called \textit{atomless} if for any $A\in\F$ with
$\PP(A)>0$, there exist $A_1,A_2\in\F$ such that
$A_1\cap A_2\ne\emp$ and $\PP(A_i)>0$.

\begin{Corollary}
\label{FR5}
Suppose that $(\Omega,\F,\PP)$ is atomless and $u$ is
law invariant. Then, for $X\in L_s^1(\DDD)$, we have
$$
u^f(X;Y)=u(\EE(X\cond Y)).
$$
\end{Corollary}

\Proof
By Corollary~A.2,
$L_s^1(\EE(\DDD\cond Y))\subseteq L_s^1(\DDD)$.
Furthermore, it is clear from~\eqref{crm12} that
$\PP\in\DDD_\mu$ for any~$\mu$, so that, by Theorem~A.1,
$\PP\in\DDD$ and $L^1\subseteq L_s^1(\DDD)$.
Now, the result follows from Theorem~\ref{FR3}.\End

\begin{Example}\rm
\label{FR6}
Let $u$ be a law invariant coherent utility
that is finite on Gaussian random variables.
Let $(X,Y)=(X,Y^1,\dots,Y^M)$ have a jointly Gaussian
distribution. Denote $\wl X=X-\EE X$, $\wl Y=Y-\EE Y$.
We can represent~$X$ as
$$
X=\EE X+\sum_{m=1}^M h^m\wl Y^m+\xi,
$$
where $\EE\xi=0$ and $\xi$ is independent of~$Y$.
Then
$$
u^f(X;Y)
=\EE X-\ga\Bigl(\EE\Bigl(
\sum_{m=1}^M h^m\wl Y^m\Bigr)^2\Bigr)^{1/2}
=\EE X-\ga\bigl|\pr_{\Lin(\wl Y^1,\dots,\wl Y^M)}(\wl X)\bigr|,
$$
where $\Lin(\wl Y^1,\dots,\wl Y^M)$ denotes the linear
space spanned by $\wl Y^1,\dots,\wl Y^M$, $\pr$ denotes the
projection, and $\ga$ is provided by Example~\ref{CRM11}~(i).

Let us give an expression for $u^f(X;Y)$ in the
matrix form. Denote $a^m=\cov(X,Y^m)$ and let $C$ be
the covariance matrix of~$Y$.
Let $L$ be the image of~$\R^M$ under the map
$x\mapsto Cx$. Note that $a\in L$.
The inverse $C^{-1}:L\to L$ is correctly defined.
The vector $h$ is found from the condition
$$
\cov(X^m-\lb h,Y\rb,Y^m)=a^m-(Ch)^m=0,\quad m=1,\dots,M,
$$
which shows that $h=C^{-1}a$. We have
$$
\bigl|\pr_{\Lin(\wl Y^1,\dots,\wl Y^M)}(\wl X)\bigr|^2
=\EE\lb h,\wl Y\rb^2
=\lb h,Ch\rb
=\lb C^{-1}a,a\rb.
$$
As a result,
$$
u^f(X;Y)=\EE X-\ga\lb C^{-1}a,a\rb^{1/2}.
$$
In particular, if $Y$ is one-dimensional, then
$$
u^f(X;Y)=\EE X-\ga\,\frac{|\cov(X,Y)|}{(\var Y)^{1/2}}.
$$

\vspace{-8mm}
\hfill$\Box$
\end{Example}

\vspace{5mm}
From the financial point of view, $\rho^f(X;Y)$ means
the risk of~$X$ in view of the uncertainty contained
in~$Y$.
So, we could expect that passing from $Y$ to a
compound random vector $(Y,Y')$ would increase~$\rho^f$.
The theorem below states that this is indeed true
for law invariant risk measures.

\begin{Theorem}
\label{FR7}
Suppose that $(\Omega,\F,\PP)$ is atomless and $u$ is
law invariant. Then, for any random variable~$X$ and
any random vectors $Y,Y'$, we have
$$
u^f(X;Y,Y')\le u^f(X;Y).
$$
\end{Theorem}

\Proof
By Theorem~A.1,
$\DDD=\bigcup_{\mu\in\mathfrak M}\DDD_\mu$ with some set
$\mathfrak M$ of probability measures on $(0,1]$.
It is seen from~\eqref{crm12} and the Jensen inequality
that
\begin{align*}
\EE(\DDD_\mu\cond Y)
&=\{Z\in L^0:Z\text{ is } Y\text{-measurable, }
Z\ge0,\;\EE Z=1,\\
&\hspace*{7mm}\text{\rm and }
\EE(Z-x)^+\le\Phi_\mu(x)\;\forall x\in\R_+\},
\end{align*}
and the similar representation is true for
$\EE(\DDD_\mu\cond Y,Y')$.
Now, it is clear from the Jensen inequality that
$\EE(\DDD_\mu\cond Y)\subseteq\EE(\DDD_\mu\cond Y,Y')$,
so that $\EE(\DDD\cond Y)\subseteq\EE(\DDD\cond Y,Y')$,
and the result follows.\End

\skm
\Remark
Coherent risk measures $\rho$ with the property
$\rho(\EE(X\mid\G))\le\rho(X)$ for any~$X$ and any
sub-$\si$-field $\G$ of~$\F$ are called
\textit{dilatation monotonous} (this property was
introduced by J.~Leitner~\cite{L04}).
Thus, Theorems~\ref{FR3} and~\ref{FR7} show that
on an atomless space any law invariant risk measure is
dilatation monotonous.\End

\skm
The example below shows that the condition of law invariance
is important in Theorem~\ref{FR7}.

\begin{Example}\rm
\label{FR8}
Let $\DDD$ consist of a unique measure~$\QQ$.
Take $Y=0$ and let $Y'$ be such that $\si(Y')=\F$.
Then $u^f(X;Y)=\EE X$, while
$u^f(X;Y,Y')=u(X)=\EE_\QQ X$.
Clearly, the inequality $\EE_\QQ X\le\EE X$ might
be violated.\End
\end{Example}

\textbf{3. Factor model.}
The question that immediately arises in connection with
the factor risks is: How close is $u^f(X;Y)$ to $u(X)$?
Below we provide a sufficient condition for the
closeness between $u^f(X;Y)$ and $u(X)$. This will be
done within the framework of the factor model, which is
very popular in statistics.

Let $u$ be a law invariant coherent utility
that is finite on Gaussian random variables and
$F=(F^1,\dots,F^M)$ be a random vector whose components
belong to $L_w^1(\DDD)$.
We will assume that $u(\lb b,F\rb)<0$
for any $b\in\R^M\setminus\{0\}$. Let
$$
X_n=B_nF+\xi_n,\quad n\in\N,
$$
where $B_n$ is $n\times M$-matrix
and $\xi_n=(\xi_n^1,\dots,\xi_n^n)$, where
$\xi_n^1,\dots,\xi_n^n,F$ are independent and
$\xi_n^i$ is Gaussian with mean~0 and variance
$(\si_n^i)^2$. We will assume that there exists a
sequence $(a_n)$ such that $a_n\to\infty$ and
\begin{align*}
&a_n^{-1}\sum_{k=1}^n B_n^{km}\xra[n\to\infty]{}b^m,
\quad m=1,\dots,M,\\
&a_n^{-2}\sum_{k=1}^n(\si_n^k)^2\xra[n\to\infty]{}0
\end{align*}
with $b=(b^1,\dots,b^M)\ne0$.

\begin{Theorem}
\label{FR9}
We have
$$
\frac{u^f\bigl(\sum_{k=1}^n X_n^k;F\bigr)}%
{u\bigl(\sum_{k=1}^n X_n^k\bigr)}\xra[n\to\infty]{}1.
$$
\end{Theorem}

\Proof
We have
$$
a_n^{-1}u\Bigl(\sum_{k=1}^n X_n^k\Bigr)
=u(\lb b_n,F\rb+\eta_n),
$$
where $b_n^m=a_n^{-1}\sum_{k=1}^n B_n^{km}$,
$m=1,\dots,M$, and $\eta_n$ is Gaussian
with mean~0 and variance
$\si_n^2=a_n^{-2}\sum_{k=1}^n(\si_n^k)^2$.
Moreover, $\eta_n$ is independent of~$F$.
Let $\eta$ be a Gaussian random variable with mean~0
and variance~1 that is independent of~$F$.
In view of the law invariance of~$u$,
$$
u(\lb b_n,F\rb+\eta_n)=u(\lb b_n,F\rb+\si_n\eta).
$$
Consider the set $G=\cl\{\EE_\QQ(F,\eta):\QQ\in\DDD\}$,
where ``$\cl$'' denotes the closure and $\DDD$ is the
determining set of~$u$.
In view of the inclusions $F^m\in L_w^1(\DDD)$,
$\eta\in L_w^1(\DDD)$, the set $G$ is a convex compact
in $\R^{M+1}$ (in the terminology of~\cite{C061},
$G$ is the \textit{generator} of $(F,\eta)$ and~$u$).
Then
\begin{align*}
&u(\lb b_n,F\rb+\si_n\eta)
=\inf_{\QQ\in\DDD}\EE_\QQ(\lb b_n,F\rb+\si_n\eta)
=\inf_{x\in G}\lb(b_n,\si_n),x\rb\\
&\xra[n\to\infty]{}\inf_{x\in G}\lb(b,0),x\rb
=\inf_{\QQ\in\DDD}\EE_\QQ\lb b,F\rb
=u(\lb b,F\rb).
\end{align*}
In a similar way we prove that
$$
a_n^{-1}u^f\Bigl(\sum_{k=1}^n X_n^k;F\Bigr)
=u(\lb b_n,F\rb)
\xra[n\to\infty]{}u(\lb b,F\rb).
$$
To complete the proof, it is sufficient to note that
$u(\lb b,F\rb)<0$.\End

\skb
\textbf{4. Multifactor risk.}
The previous theorem shows that in order to assess the
risk $\rho(X)$, one can take the main factors
$Y^1,\dots,Y^M$ driving risk, and then
$\rho(X)\approx\rho^f(X;Y^1,\dots,Y^M)$.
The question arises: What is the relationship between
$\rho^f(X;Y^1,\dots,Y^M)$ and
$\rho^f(X;Y^1)+\dots+\rho^f(X;Y^M)$?
First, we provide a positive statement.

\begin{Proposition}
\label{FR10}
Assume that $Y^1,\dots,Y^M$ are independent and
$X=\sum_{m=1}^MX^m$, where $X^m$ is $Y^m$-measurable,
$X^m\in L_s^1(\DDD)\cap L_s^1(\EE(\DDD\cond Y))\cap L^1$,
and $\EE X^m=0$. Then
\begin{equation}
\label{fr2}
\rho^f(X;Y^1,\dots,Y^M)\le\sum_{m=1}^M\rho^f(X;Y^m).
\end{equation}
\end{Proposition}

\Proof
We have
\begin{align*}
&\rho^f(X;Y^1,\dots,Y^M)
=\rho(\EE(X\cond Y^1,\dots,Y^M))
=\rho(X),\\
&\rho^f(X;Y^m)
=\rho(\EE(X\mid Y^m))
=\rho(X^m),
\end{align*}
and the result follows from the subadditivity property
of~$\rho$.\End

\skm
The conditions of the proposition are unrealistic
because different factors are correlated.
This might lead to violation of~\eqref{fr2} as shown
by the example below.

\begin{Example}\rm
\label{FR11}
Let $(Y^1,Y^2)$ be a Gaussian random vector with mean~0
and covariance matrix
$$
\begin{pmatrix}
1&1-\eps\\
1-\eps&1
\end{pmatrix}.
$$
Take $X=Y^1-Y^2$.
Let $\rho$ be a law invariant coherent risk measure
that is finite on Gaussian random variables.
Then
$$
\rho^f(X;Y^1,Y^2)
=\rho(X)
=\ga(\var(Y^1-Y^2))^{1/2}
=\ga\sqrt{2\eps},
$$
where $\ga$ is provided by Example~\ref{CRM11}~(i).
On the other hand, by Example~\ref{FR6},
$$
\rho^f(X;Y^n)=
\ga\,\frac{|\cov(X,Y^n)|}{(\var Y^n)^{1/2}}
=\ga\eps,\quad n=1,2.
$$
If $\eps$ is small enough, then the inequality
$\ga\sqrt{2\eps}\le2\ga\eps$ is violated.\End
\end{Example}

\begin{figure}[!h]
\begin{picture}(150,56)(-70,-32.5)
\put(-15,-10){\includegraphics{factor.4}}
\put(31,9){\small $Y^1$}
\put(31,-11){\small $Y^2$}
\put(-1.5,21){\small $X$}
\put(-30,-25){\parbox{80mm}{\small\textbf{Figure~6.}
In Example~\ref{FR11},
$\rho^f(X;Y^1;Y^2)/\ga$ is the length of~$X$,
while $\rho^f(X;Y^n)/\ga$ is the length of the projection
of~$X$ on~$Y^n$.}}
\end{picture}
\end{figure}

The effect described in this example has the following
financial background.
Suppose that we have several correlated factors
$Y^1,\dots,Y^M$ and a portfolio consisting of~$M$
parts, i.e. $X=\sum_m X^m$, where the risk of the $m$-th
part is driven mainly by the $m$-th factor.
Then $X^m$ is correlated with $Y^k$ for $m\ne k$
simply because $Y^m$ and $Y^k$ are correlated.
Thus, when summing up $\rho^f(X;Y^m)$ over~$m$,
we are calculating the factor loading of
$Y^m$ in $X^m$ once through the $m$-th factor risk and
then several times more through the other correlated
factor risks.
This might lead to a significant increment as well as a significant
reduction of the estimated risk (which was described by
Example~\ref{FR11}).
In this situation, the right way to estimate $\rho(X)$
is to take $\rho^f(X^1;Y^1)+\dots+\rho^f(X^M;Y^M)$.
Indeed, if each $X^m$ corresponds to a big portfolio, then
above results tell us that $\rho(X^m;Y^m)\approx\rho(X^m)$,
and then
$$
\rho(X)\le\rho\Bigl(\sum_{m=1}^M X^m\Bigr)
\approx\sum_{m=1}^M\rho^f(X^m;Y^m).
$$
Another pleasant feature of this technique is that we
estimate the $m$-th factor risk only for a part of the
portfolio rather than the whole portfolio, which
accelerates the computation speed.

However, it might happen that we cannot split a portfolio
in several groups such that the risk of each group is
affected by one factor only (for example, the risk of
credit derivatives is affected by the whole yield curve).
But then we might combine the one-factor and the multifactor
techniques as follows.
Suppose that we can split the factors into several groups
(thus we again have $Y^1,\dots,Y^M$, but now these
random variables are multidimensional) and split the portfolio
into~$M$ parts, i.e. $X=\sum_m X^m$, so that the risk of
the $m$-th part is driven mainly by the $m$-th group of
factors. Then we can assess the risk of the portfolio as
$\rho^f(X^1;Y^1)+\dots+\rho^f(X^M;Y^M)$.

\skm
\textbf{5. Empirical estimation.}
In view of Theorem~\ref{FR3}, the empirical estimation
of $\rho^f(X;Y)$ reduces to finding the function
$f(y)=\EE(X\cond Y=y)$ and then applying the procedures
described at the end of Section~\ref{CRM} to $x_{kl}:=f(y_{kl})$.
However, for factor risks we can use one more
convenient method of the choice of data.

Suppose that $Y$ is one-dimensional and let $\si$ be
the current volatility of~$Y$. It might be estimated
through one of numerous well-known methods; in particular,
$\si$ might be the implied volatility.
It is a widely accepted idea that volatility serves
as the speed of growth of the inner time
(known also as the business or operational time)
for~$Y$. In other words, if the current volatility
is~$\si$, then $Y$ is currently oscillating at the
speed~$\si^2$. Following this idea, we can take as the
data for~$Y$ the values $y_1,\dots,y_T$, where
$y_t$ is the increment of~$Y$ over the time interval
$[-t\si^2\De,-(t-1)\si^2\De]$ and $\De$ is the unit time
period.
Here it is reasonable to take standardized time series
for~$Y$ rather than the ordinary one:
we calculate empirically the integrated volatility
and take its inverse as the time change to obtain
standardized time series from the ordinary one.
This approach enables one
\begin{mitemize}
\item to use large data sets;
\item to capture volatility predictions immediately.
\end{mitemize}
If $\De$ is one day (which is a typical choice), then
$\si^2\De$ would be a non-integer number of days,
which is not very convenient. This can be overcome as follows.
We choose a large number~$n\in\N$ and split the time axis
into small intervals of length $n^{-1}\De$.
Then we approximate $\si^2$ by a rational number $m/n$
and generate each $y_t$ as a sum of increments of~$Y$
over $m$ randomly chosen small intervals.
In other words, this is a combination of the time change
procedure with the bootstrap technique.

Instead of the time change procedure described above,
one can also use the classical scaling procedure, which
is at the basis of the \textit{filtered historical
simulation} (see~\cite[Sect.~5.6]{C03}).
Namely, instead of altering the time step for~$y_t$,
we keep the same time step~$\De$, but multiply each~$y_t$
by the current volatility~$\si$.
As above, it is reasonable to take standardized time series
for~$Y$ rather than the ordinary one:
we divide each observed increment of~$Y$ by its
volatility estimated through one of standard techniques
(for example, GARCH).

Let us remark that both the time change and
the scaling work
only for one-dimensional $Y$s because if $Y$ is multidimensional,
its different components have different volatilities.
This is a big advantage of one-factor risks.

%======================================================
\section{Portfolio Optimization}
\label{PO}

\textbf{1. Problem.}
Let $(\Omega,\F,\PP)$ be a probability space and
$u^1,\dots,u^M$ be coherent utilities with the
determining sets $\DDD^1,\dots,\DDD^M$.
A particular example we have in mind is
$u^m=u^f(\,\cdot\,;Y^m)$, where $Y^1,\dots,Y^M$ are the
main factors driving the risk of a portfolio.
Let $X^1,\dots,X^d\in\bigcap_m L_w^1(\DDD^m)$ be the
P\&Ls produced by traded assets over the unit time period,
so that the space of possible P\&Ls attained by various
investment strategies is $\{\lb h,X\rb:h\in\R^d\}$, where
$X=(X^1,\dots,X^d)$.
We will assume that $u^m(\lb h,X\rb)<0$ for any
$h\in\R^d\setminus\{0\}$.
This condition means that the risk of any trade is
strictly positive and is known as the \textit{No Good Deals}
condition (see~\cite{C061}).
Let $E=(E^1,\dots,E^d)$ be the vector of rewards
for $X^1,\dots,X^d$. One might think of~$E$ as the vector
of expectations $(\EE X^1,\dots,\EE X^d)$. However, this
is not the only interpretation of~$E$. In general, we
mean by~$E$ the vector of subjective assessments by
some agent of the profitability of the assets
$X^1,\dots,X^d$. In this case $E$ need not be related
to $\EE X$, and $\EE X$ can be equal to zero.

We will consider the Markowitz-type optimization problem,
risk being measured not as variance, but rather
as the vector of risks:
\begin{equation}
\label{po1}
\begin{cases}
\lb h,E\rb\longrightarrow\max,\\
h\in\R^d,\\
\rho^m(\lb h,X\rb)\le c^m,\;m=1,\dots,M,
\end{cases}
\end{equation}
where $c^1,\dots,c^M\in(0,\infty)$ are fixed risk limits.

\skb
\textbf{2. Geometric solution.}
The paper~\cite[Subsect.~2.2]{C062} contains a geometric
solution of this problem with $M=1$.
Here we will present a similar geometric solution
for an arbitrary~$M$.
Let us introduce the notation
\begin{align*}
G^m&=\cl\{\EE_\QQ X:\QQ\in\DDD^m\},\quad m=1,\dots,M,\\
G&=\conv\{G^m/c^m:m=1,\dots,M\}.
\end{align*}
Note that $G^m,G$ are convex compacts in $\R^d$.
According to the terminology of~\cite{C061}, $G^m$
is the \textit{generator} of $X$ and $u^m$.
The role of this set is seen from the equality
\begin{equation}
\label{po2}
u^m(\lb h,X\rb)=\inf_{x\in G^m}\lb h,x\rb,\quad m=1,\dots,M.
\end{equation}
The right-hand side is the classical object of
convex analysis termed the \textit{support function}
of the set~$G^m$.
The notion of the generator was found to be very
convenient for the geometric solutions of various
problems like capital allocation, portfolio optimization,
pricing, and equilibrium (see~\cite{C061}, \cite{C062}).
As $u^m(\lb h,X\rb)<0$ for any $h\in\R^d\setminus\{0\}$,
the point~0 belongs to the interior of~$G$.
Let $T$ be the intersection of the ray $(E,0)$
with the border of~$G$. Denote by $N$ the set of inner
normals to~$G$ at the point~$T$ (typically, $N$ is a ray).

\begin{figure}[!h]
\begin{picture}(150,92.5)(-39,-30)
\put(0,0){\rotatebox{345}{\includegraphics{factor.5}}}
\put(22,26.5){\small $G$}
\put(56,32.8){\small $\frac{G^1}{c^1}$}
\put(48,4.5){\small $\frac{G^M}{c^M}$}
\put(66.5,38.5){\small $G^1$}
\put(53.7,-6){\small $G^M$}
\put(27.5,16.5){$h_*$}
\put(25.5,2){$T$}
\put(41.2,21){$0$}
\put(43.5,25){$E$}
\put(5,-25){\parbox{70mm}{\small\textbf{Figure~7.}
Solution of the optimization problem.
Here $h_*$ denotes the optimal $h$.}}
\end{picture}
\end{figure}

\begin{Theorem}
\label{PO1}
The set of solutions of problem~\eqref{po1} is
$\{h\in N:\lb h,T\rb=-1\}$ and
the maximal $\lb h,E\rb$ is $|E|/|T|$.
\end{Theorem}

\Remark
Note that $N$ is a non-empty cone, so that
$\{h\in N:\lb h,T\rb=-1\}\ne\emp$.
Furthermore, $0<|E|/|T|<\infty$.\End

\skm
\texttt{Proof of Theorem~\ref{PO1}.}
Using~\eqref{po2}, we can write
$$
\eqref{po1}
\;\Lea\;
\begin{cases}
\lb h,E\rb\longrightarrow\max,\\
h\in\R^d,\\
\inf\limits_{x\in G^m}\lb h,x\rb\ge-c^m
\end{cases}
\hspace{-3mm}\Lea\;
\begin{cases}
\lb h,E\rb\longrightarrow\max,\\
h\in\R^d,\\
\inf\limits_{x\in G}\lb h,x\rb\ge-1.
\end{cases}
$$

Denote $\{h\in N:\lb h,T\rb=-1\}$ by $H_*$.
For any $h\in H_*$, we have
$$
\inf_{x\in G}\lb h,x\rb=\lb h,T\rb=-1
$$
and
\begin{equation}
\label{po3}
\lb h,E\rb
=-\frac{|E|}{|T|}\lb h,T\rb
=\frac{|E|}{|T|}.
\end{equation}

If $h\notin H_*$ and $\inf_{x\in G}\lb h,x\rb\ge-1$,
then $\inf_{x\in G}\lb h,x\rb<\lb h,T\rb$, so that
$\lb h,T\rb>-1$, and, due to~\eqref{po3},
$\lb h,E\rb<|E|/|T|$.\End

\skm
\Remark
We can also provide a geometric solution of~\eqref{po1}
under portfolio constraints of the type $h\in H$, where
$H$ is a convex cone, and the ambiguity of the reward
vector~$E$.
This is done by transforming~\eqref{po1} into the
problem with one constraint as described above and
then applying the result of~\cite[Subsect.~2.2]{C062},
where the cone constraints and the ambiguity were taken
into account.\End

\begin{Example}\rm
\label{PO2}
Let $X^1,\dots,X^d\in L^1(\DDD_\mu)$
and $Y$ be an $M$-dimensional random vector.
Then the generator of $u_\mu^f(\,\cdot\,;Y)$ and~$X$
has the form
\begin{align*}
\cl\{\EE_\QQ X:\QQ\in\EE(\DDD_\mu\cond Y)\}
&=\cl\Bigl\{\int_{\R^M}f(y)Z(y)\QQ(dy):
Z\ge0,\,\int_{\R^M}Z(y)\QQ(dy)=1,\\
&\hspace*{11mm}\text{and }\int_{\R^M}(Z(y)-x)^+\QQ(dy)
\le\Phi_\mu(x)\;\forall x\in\R_+\Bigr\},
\end{align*}
where $f(y)=\EE(X\cond Y=y)$, $\QQ=\Law Y$, and
$\Phi_\mu$ is given by~\eqref{crm13}.
In order to prove this equality, denote its left-hand
side by~$G$ and its right-hand side by~$G'$.
Due to~\eqref{crm12},
$$
\inf_{x\in G'}\lb h,x\rb=u'_\mu(\lb h,f\rb),
\quad h\in\R^d,
$$
where $u'_\mu$ is minus the Weighted V@R with the
weighting measure~$\mu$ on the probability space
$(\R^M,\B(\R^M),\QQ)$.
As $u_\mu$ and $u'_\mu$ depend only on the distribution
of a random variable (this is seen from~\eqref{crm6},
\eqref{crm10}) and $\Law_\QQ\lb h,f\rb=\Law\lb h,f(Y)\rb$,
we get
\begin{align*}
u'_\mu(\lb h,f\rb)
&=u_\mu(\lb h,f(Y)\rb)
=u_\mu(\EE(\lb h,X\rb\cond Y))
=u^f(X;Y)\\[1mm]
&=\inf_{\QQ\in\EE(\DDD_\mu\cond Y)}\EE_\QQ\lb h,X\rb
=\inf_{x\in G}\lb h,x\rb,\quad h\in\R^d.
\end{align*}
Thus, the support functions of $G$ and $G'$ coincide.
Furthermore, both $G$ and $G'$ are convex and closed.
As a result, $G=G'$.\End
\end{Example}

\textbf{3. Practical aspects.}
The geometric solution presented above provides a nice
theoretical insight into the form of the optimal
portfolio. It can be used if we have a model for the
joint distribution of $X^1,\dots,X^d$. However, if we
do not have such a model, but rather want to
approach~\eqref{po1} empirically, then, instead of the
geometric solution, the following straightforward
procedure can be employed. First of all,
$$
\eqref{po1}
\;\Lea\;
\begin{cases}
\lb h,E\rb\longrightarrow\max,\\
h\in\R^d,\\
\ds\min_{m=1,\dots,M}\frac{u^m(\lb h,X\rb)}{c^m}\ge-1.
\end{cases}
$$
Due to the scaling property, the solution to this problem
coincides up to multiplication by a positive (easily
computable) constant with the solution of the problem
$$
\begin{cases}
\ds\min_{m=1,\dots,M}\frac{u^m(\lb h,X\rb)}{c^m}
\longrightarrow\max,\\
h\in\R^d,\\
\lb h,E\rb=1.
\end{cases}
$$
This is a problem of maximizing a convex functional
over an affine space.
A typical example we have in mind is
$u^m=u^f(\,\cdot\,;Y^m)$, where $u$ is one of $u_\mu$,
$u_{\al,\be}$, or $u_\al$. In this case the values
$u^m(\lb h,X\rb)$ can easily be estimated empirically
through the procedures described at the end of
Section~\ref{CRM}.

If $E$ is the vector of expected profits
$(\EE X^1,\dots,\EE X^d)$, then its empirical
estimation is known to be an extremely unpleasant
problem (see the discussion in~\cite{B95} and
the 20's example in~\cite{JPR05}), unlike the estimation
of the volatility-type quantities $u^m$.
The reason is that this vector is very close to zero,
and therefore, its direction (which is in fact the input
that we need) depends on the data in a very unstable way.
One of possible ways to overcome this problem is to
use theoretical estimates of~$E$ rather than
empirical ones.
For example, Sharpe's SML relation (\cite{S64}) implies that
$$
E^i=\beta^i\times const,\quad i=1,\dots,d
$$
(recall that $X^i$ is the \textit{discounted}
P\&L produced by the $i$-th asset).
Thus, the direction of~$E$ (and this is what we need)
coincides with that of $(\beta^1,\dots,\beta^d)$.

%=======================================================
\section{Risk Contribution}
\label{RC}

\textbf{1. Extreme measures.}
Let $(\Omega,\F,\PP)$ be a probability space
and $u$ be a coherent utility with the
determining set~$\DDD$.
Let $W$ be a random variable meaning the P\&L produced
by some portfolio over the unit time period.

The following definition was introduced in~\cite{C061}.

\begin{Definition}\rm
\label{RC1}
A measure $\QQ\in\DDD$ is an \textit{extreme measure} for~$W$ if
$\EE_\QQ W=u(W)\in(-\infty,\infty)$.

The set of extreme measures will be denoted by $\EX_\DDD(W)$.
\end{Definition}

\begin{Proposition}
\label{RC2}
If the determining set $\DDD$ is $L^1$-closed and
uniformly integrable, while $W\in L_s^1(\DDD)$,
then $\EX_\DDD(W)\ne\emp$.
\end{Proposition}

\Proof
It is clear that $u(X)\in(-\infty,\infty)$.
Find a sequence $Z_n\in\DDD$ such that
$\EE Z_n X\to u(X)$.
By the Dunford--Pettis criterion, $\DDD$ is compact
with respect to the weak topology $\si(L^1,L^\infty)$.
Therefore, the sequence $(Z_n)$ has a weak limit point
$Z_\infty\in\DDD$.
Clearly, the map $\DDD\ni Z\mapsto\EE ZX$ is
weakly continuous.
Hence, $\EE Z_\infty X=u(X)$, which means that
$Z_\infty\in\EX_\DDD(X)$.\End

\skm
\Remark
It is seen from~\eqref{crm12} that the determining set
of Weighted V@R is $L^1$-closed and uniformly integrable
(note that $\Phi_\mu(x)\to0$ as $x\to\infty$).\End

\begin{Example}\rm
\label{RC3}
\textbf{(i)} If $\la\in(0,1]$ and $W\in L^1$, then
\begin{equation}
\begin{split}
\label{rc1}
\EX_{\DDD_\la}(W)
&=\bigl\{Z:Z\ge0,\,\EE Z=1,\,Z=\la^{-1}\text{ a.s. on }
\{W<q_\la(W)\},\\
&\hspace*{7mm}\text{and }Z=0\text{ a.s. on }
\{W>q_\la(W)\}\bigr\},
\end{split}
\end{equation}
where $\DDD_\la$ is given by~\eqref{crm3}.
Indeed, if $Z$ belongs to the right-hand side of~\eqref{rc1},
then, for any $Z'\in\DDD_\la$, we have
\begin{align*}
\EE W Z'-\EE W Z
&=\EE(W-q_\la(W))Z'-\EE(W-q_\la(W))Z\\
&=\EE[(W-q_\la(W))(Z'-\la^{-1})I(W-q_\la(W)<0)\\
&\hspace*{8mm}+(W-q_\la(W))Z'I(W-q_\la(W)>0)].
\end{align*}
It is seen that this quantity is positive and equals
zero if and only if $Z'$ belongs to the right-hand side
of~\eqref{rc1}.

\textbf{(ii)} Let $\Omega=\{1,\dots,T\}$ and
$W(t)=w_t$.
Assume that $w_1<\dots<w_T$.
Set $z_t=\sum_{i=1}^t\PP\{i\}$.
It is seen from~(i) that, for any $\la\in(0,1]$,
$\EX_{\DDD_\la}(W)$ consists of the unique measure
$\QQ_\la(W)$ having the form
$$
\QQ_\la(W)\{i\}=\begin{cases}
\la^{-1}\PP\{i\},&i\le n-1,\\
1-\la^{-1}z_{n-1},&i=n,\\
0,&i>n,
\end{cases}
$$
where $n$ is such that $z_{n-1}<\la\le z_n$.
This can be rewritten as
$$
\QQ_\la(W)\{i\}=\la^{-1}\int_0^\la I(z_{i-1}<x\le z_i)dx,
\quad i=1,\dots,T.
$$
It follows from~\cite[Prop.~6.2]{C05e} that
$\EX_{\DDD_\mu}(W)$ consists of the unique measure
$\QQ_\mu(W)=\int_{(0,1]}\QQ_\la(W)\mu(d\la)$. We have
\begin{align*}
\QQ_\mu(W)\{i\}
&=\int_{(0,1]}\int_0^\la\la^{-1}I(z_{i-1}<x\le z_i)
dx\mu(d\la)\\[1mm]
&=\int_0^1\int_{[x,1]}\la^{-1}I(z_{i-1}<x\le z_i)
\mu(d\la)dx\\[1mm]
&=\int_0^1 I(z_{i-1}<x\le z_n)\psi_\mu(x)dx\\[1mm]
&=\int_{z_{n-1}}^{z_n}\psi_\mu(x)dx,
\quad i=1,\dots,T,
\end{align*}
where $\psi_\mu$ is given by~\eqref{crm7}.

\textbf{(iii)} If $W\in L^1(\DDD_\mu)$ has a continuous
distribution, then $\EX_{\DDD_\mu}(W)$ consists of the
unique measure $\QQ_\mu(W)=\psi_\mu(F(W))\PP$, where
$F$ is the distribution function of~$W$ (for the proof,
see~\cite[Sect.~6]{C05e}). The measure $\QQ_\mu$ has
already appeared in~\eqref{crm11}.

\textbf{(iv)} Suppose that $\al,\be\in\N$ and
$W\in L^1(\DDD_\mu)$ has a continuous distribution.
For Beta V@R, the measure $\QQ_\mu(W)$ of the previous
example gets a more concrete form
$\QQ_{\al,\be}(W)
=\psi_{\al,\be}(F(W))\PP
=\phi(W)\PP$, where
$\psi_{\al,\be}$ is provided by~\eqref{crm19}.

Let us clarify the meaning of $\phi$.
Let $W_1,\dots,W_\al$ be independent copies of~$W$,
$W_{(1)},\dots,W_{(\al)}$ be the corresponding order
statistics, and $\xi$ be an independent uniformly distributed
on $\{1,\dots,\be\}$ random variable.
According to the reasoning of Example~\ref{CRM9},
the distribution function of $W_{(\xi)}$ is
$\Psi_{\al,\be}\circ F$, where
$$
\Psi_{\al,\be}(x)=\int_0^x\psi_{\al,\be}(y)dy,
\quad x\in[0,1].
$$
Consequently,
$$
\frac{d\Law W_{(\xi)}}{d\Law W}\,(x)
=\frac{d\Psi_{\al,\be}(F(x))}{dF(x)}
=\psi_{\al,\be}(F(x))
=\phi(x).
$$

\textbf{(v)} For Alpha V@R with $\al\in\N$,
we have $\psi_\al(x)=\al(1-x)^{\al-1}$ and
$$
\phi(x)
=\al(1-F(x))^{\al-1}
=\frac{d\Law\min\{W_1,\dots,W_\al\}}{d\Law W}.
$$

\vspace{-8mm}
\hfill$\Box$
\end{Example}

\vspace{5mm}
If an agent is using the classical expected utility
$\EE U(X)$ to assess the quality of his/her position,
then there exists his/her ``personal'' measure with which
he/she assesses the quality of any possible trade.
This measure is given by $\QQ=cU'(W_1)\PP$, where
$W_1$ is the agent's wealth at the terminal date
and $c$ is the normalizing constant.
The role of this measure is seen from the equality
$$
\lim_{\eps\da0}\eps^{-1}(\EE U(W_1+\eps X)-\EE U(W_1))
=\EE XU'(W_1)
=c^{-1}\EE_\QQ X.
$$
If $X$ is the P\&L produced by some trade and
$X$ is small as compared to~$W$, then $X$ is profitable
for the agent if and only if $\EE_\QQ X>0$.
The extreme measure is the coherent substitute for this
agent-specific measure as seen from
Theorem~\ref{RC5} stated below.

\skb
\textbf{2. Risk contribution.}
The following definition was introduced in~\cite{C061}.

\begin{Definition}\rm
\label{RC4}
The \textit{utility contribution} is defined as
$$
u^c(X;W)=\inf_{\QQ\in\EX_\DDD(W)}\EE_\QQ X,\quad X\in L^0.
$$
The \textit{risk contribution} is defined as
$\rho^c(X;W)=-u^c(X;W)$.
\end{Definition}

If $\EX_\DDD(W)$ is non-empty, then $u^c(\,\cdot\,;W)$
is a coherent utility.

\begin{Theorem}
\label{RC5}
If $\DDD$ is $L^1$-closed and uniformly integrable
and $X,W\in L_s^1(\DDD)$, then
$$
u^c(X;W)=\lim_{\eps\da0}\eps^{-1}(u(W+\eps X)-u(W)).
$$
\end{Theorem}

For the proof, see~\cite[Subsect.~2.5]{C061}.

\skm
For theoretical purposes, it is sometimes convenient to
use a geometric representation of~$u^c(X;W)$.
Suppose that $\DDD$ is $L^1$-closed and uniformly
integrable, while $X,W\in L_s^1(\DDD)$.
Denote by~$G$ the generator
$\cl\{\EE_\QQ(X,W):\QQ\in\DDD\}$ and set $e_1=(1,0)$,
$e_2=(0,1)$. Then
$$
\EX_\DDD(W)=\bigl\{\QQ\in\DDD:\EE_\QQ W=\min_{z\in G}
\lb e_2,z\rb\bigr\},
$$
and therefore,
\begin{equation}
\label{rc2}
u^c(X;W)=\min\bigl\{x:\bigl(x,\min_{z\in G}
\lb e_2,z\rb\bigr)\in G\bigr\}.
\end{equation}

\begin{figure}[!h]
\begin{picture}(150,59)(-55,-16.5)
\put(4.7,5){\includegraphics{factor.6}}
\put(0,0){\vector(1,0){50}}
\put(0,0){\vector(0,1){40}}
\multiput(0,5)(2,0){10}{\line(1,0){1}}
\multiput(5,0)(0,2){10}{\line(0,1){1}}
\multiput(20,0)(0,2){3}{\line(0,1){1}}
\put(48,-4){\small $x$}
\put(-3.5,38){\small $y$}
\put(-9,4){\scalebox{0.8}{$u(W)$}}
\put(1,-4){\scalebox{0.8}{$u(X)$}}
\put(14,-4){\scalebox{0.8}{$u^c(X;W)$}}
\put(23,19){\small $G$}
\put(-17,-14){\small\textbf{Figure~8.} Geometric representation
of~$u^c(X;W)$}
\end{picture}
\end{figure}

\begin{Example}\rm
\label{RC6}
\textbf{(i)} If $W$ is a constant,
then $\EX_\DDD(W)=\DDD$,
so that $u^c(X;W)=u(X)$.

\textbf{(ii)} If $X=\la W$ with $\la\in\R_+$,
then $u^c(X;W)=\la u(W)$
provided that $\EX_\DDD(W)\ne\emp$.

\textbf{(iii)} Let $u$ be law invariant and $(X,W)$ be
jointly Gaussian. We assume that $X,W\in L_s^1(\DDD)$
and that the covariance matrix~$C$ of $(X,W)$ is
non-degenerate. Set $\wl X=X-\EE X$, $\wl W=W-\EE W$.
As $(\wl X,\wl W)$ can be represented as $C^{1/2}V$,
where $V$ is a standard two-dimensional Gaussian random
vector, the generator~$\wl G$ of $(\wl X,\wl W)$ has the
form $\wl G=C^{1/2}G_V$, where $G_V$ is the generator
of~$V$. Clearly, $G_V$ is the ball of radius~$\ga$,
where $\ga$ is provided by Example~\ref{CRM11}~(i). Thus,
$$
\wl G
=\{x\in\R^2:\lb C^{-1/2}x,C^{-1/2}x\rb\le\ga^2\}
=\{x\in\R^2:\lb x,C^{-1}x\rb\le\ga^2\}.
$$
Set $e_1=(1,0)$, $e_2=(0,1)$. By~\eqref{rc2},
$u^c(\wl X,\wl W)=\lb e_1,z_*\rb$, where
$z_*=\argmin_{z\in G}\lb e_2,z\rb$.
The point~$z_*$ is found from the condition
$$
0=\frac{d}{d\eps}\Bigl|_{\eps=0}
\lb z_*+\eps e_1,C^{-1}(z_*+\eps e_1)\rb
=2\lb e_1,C^{-1}z_*\rb,
$$
which shows that $C^{-1}z_*=\al e_2$ with some $\al\le0$,
i.e. $z_*=\al Ce_2$. The constant~$\al$ is found from the
condition $\lb z_*,C^{-1}z_*\rb=\ga^2$, which shows that
$\al=-\ga\lb e_2,Ce_2\rb^{-1/2}$. As a result,
$$
u^c(X;W)
=\EE X+u^c(\wl X,\wl W)\\[1mm]
=\EE X-\ga\,\frac{\lb e_1,Ce_2\rb}{\lb e_2,Ce_2\rb^{1/2}}\\[1mm]
=\EE X-\ga\,\frac{\cov(X,W)}{(\var W)^{1/2}}.
$$
Note that $u(X)=\EE X-\ga(\var X)^{1/2}$.
In particular, if $\EE X=\EE W=0$, then
\begin{equation}
\label{rc3}
\frac{u^c(X;W)}{u(X)}
=\corr(X,W)
=\frac{\text{V@R}^c(X;W)}{\text{V@R}(X)},
\end{equation}
where $\text{V@R}^c$ denotes the V@R contribution
(for the definition, see~\cite[Sect.~7]{M02}).

\textbf{(iv)} Let $\Omega=\{1,\dots,T\}$,
$X(t)=x_t$, $W(t)=w_t$.
Assume that all the values $w_t$ are different.
Let $w_{(1)},\dots,w_{(T)}$ be the values
$w_1,\dots,w_T$ in the increasing order.
Define $n(t)$ through the equality $w_{(t)}=w_{n(t)}$.
According to Example~\ref{RC3}~(ii),
$$
u_\mu^c(X;W)
=\EE_{\QQ_\mu(W)}X
=\sum_{t=1}^T x_{n(t)}\int_{z_{t-1}}^{z_t}\psi_\mu(x)dx,
$$
where $z_t=\sum_{i=1}^t\PP\{n(i)\}$ (cf.~\eqref{crm8}).

\textbf{(v)} If $W\in L^1(\DDD_\mu)$ has a continuous
distribution, then, according to Example~\ref{RC3}~(iii),
$$
u_\mu^c(X;W)=\EE_{\QQ_\mu(W)}X=\EE X\psi_\mu(F(W))
$$
(cf.~\eqref{crm11}).
Note that this value is linear in~$X$.

\textbf{(vi)} Suppose that $\al,\be\in\N$, $X,W\in L^1$,
and $W$ has a continuous distribution.
Let $(X_1,W_1),\dots,(X_\al,W_\al)$ be independent
copies of $(X,W)$ and $\xi$ be an independent uniformly
distributed on $\{1,\dots,\be\}$ random variable.
Let $W_{(1)},\dots,W_{(\al)}$ be the corresponding order
statistics.
Define random variables $n(i)$ through the equality
$W_{(i)}=W_{n(i)}$ (as $W$ has a continuous distribution,
all the values $W_1,\dots,W_\al$ are a.s. different, so
that $n(i)$ is a.s. determined uniquely).
We have (cf.~\eqref{crm18})
\begin{align*}
\EE X_{n(\xi)}
&=\EE\Bigl[\frac{1}{\be}\sum_{i=1}^\be X_{n(i)}\Bigr]\\
&=\frac{1}{\be}\sum_{i=1}^\be\sum_{j=1}^\al
\EE X_j I(n(i)=j)\\
&=\frac{1}{\be}\sum_{i=1}^\be\sum_{j=1}^\al
\EE X_j I\{i-1\text{ members of }W_1,\dots,W_\al
\text{ are smaller than }W_j\\[-2mm]
&\hspace{25mm}\text{ and }\al-i\text{ members of }W_1,\dots,W_\al
\text{ are greater than }W_j\}\\
&=\frac{1}{\be}\sum_{i=1}^\be\sum_{j=1}^\al
\EE[X_jF(W_j)^{i-1}(1-F(W_j))^{\al-i}]\\
&=\frac{\al}{\be}\sum_{i=1}^\be C_{\al-1}^{i-1}
\EE[X F(W)^{i-1}(1-F(W))^{\al-i}]\\
&=\EE X\psi_{\al,\be}(F(W))\\
&=u_{\al,\be}^c(X;W).
\end{align*}

\textbf{(vii)} Suppose that $\al\in\N$, $X,W\in L^1$,
and $W$ has a continuous distribution.
It follows from~(vi) that
$$
u_\al^c(X;W)=\EE X_{\argmin\limits_{i=1,\dots,\al}W_i},
$$
where $(X_1,W_1),\dots,(X_\al,W_\al)$ are independent
copies of $(X,W)$.\End
\end{Example}

\textbf{3. Capital allocation.}
The notion of risk contribution is closely
connected with the \textit{capital allocation} problem.
Suppose that a firm consists of several desks, i.e.
$W=\sum_{n=1}^N W^n$, where $W^n$ is the P\&L produced
by the $n$-th desk.
The following definition was introduced in~\cite{D05}.

\begin{Definition}\rm
\label{RC7}
A collection $x^1,\dots,x^N$ is a \textit{capital
allocation} between~$W^1,\dots,W^N$~if
\begin{align}
\label{rc4}
&\sum_{n=1}^N x^n=\rho\Bigl(\sum_{n=1}^N W^n\Bigr),\\
\label{rc5}
&\sum_{n=1}^N h^nx^n\le\rho\Bigl(\sum_{n=1}^N h^nX^n\Bigr)
\;\;\forall h^1,\dots,h^N\in\R_+.
\end{align}
\end{Definition}

From the financial point of view, $x^i$ means the contribution
of the $i$-th component to the total risk of the firm, or,
equivalently, the capital that should be allocated to this
component.
In order to illustrate the meaning of~\eqref{rc5},
consider the example $h^n=I(n\in J)$, where
$J$ is a subset of $\{1,\dots,N\}$.
Then~\eqref{rc5} means that the capital
allocated to a part of the firm does not exceed the risk
carried by that part.

The following theorem was established
in~\cite[Subsect.~2.4]{C061}.

\begin{Theorem}
\label{RC8}
Suppose that $\DDD$ is $L^1$-closed and uniformly
integrable, while $W^1,\dots,W^N\in L_s^1(\DDD)$.
Then the set of solutions of the capital allocation
problem has the form
$\{-\EE_\QQ(W^1,\dots,W^N):\QQ\in\EX_\DDD(W)\}$.
\end{Theorem}

If $\EX_\DDD(W)$ consists of a unique measure~$\QQ$,
then the solution of the capital allocation problem is
unique and has the form $(\rho^c(W^1;W),\dots,\rho^c(W^N;W))$.
In particular, in this case
\begin{equation}
\label{rc6}
\rho(W)
=-\EE_\QQ W
=-\sum_{n=1}^N\EE_\QQ W^n
=\sum_{n=1}^N\rho^c(W^n;W).
\end{equation}

\skb
\textbf{4. Tail correlation.}
It follows from the inclusion $\EX_\DDD(W)\subseteq\DDD$
that $u^c(X;W)\ge u(X)$.
Typically, for a random variable~$X$ meaning the P\&L
of some transaction, we have $u(X)<0$ (i.e. its risk is
strictly positive), so that
$$
\kappa(X;W):=\frac{u^c(X;W)}{u(X)}\le1.
$$
The coefficient $\kappa$ is a good measure for the
tail correlation between~$X$ and~$W$ (see~\eqref{rc3}).

Let us recall that the standard \textit{tail correlation}
coefficient between $X$ and $W$ (see~\cite[Sect.~5.2.3]{EFM05})
is defined as $\lim_{\la\da0}c_\la(X;W)$, where
$$
c_\la(X;W)=\frac{\PP(X\le q_\la(X),\,W\le q_\la(W))}%
{\PP(X\le q_\la(X))}
=\frac{\EE I(X\le q_\la(X))I(W\le q_\la(W))}%
{\EE I(X\le q_\la(X))I(X\le q_\la(X))}.
$$
At the same time, $\kappa(X;W)$ corresponding to
$u=u_\la$ has the form
$$
\kappa_\la(X;W)=\frac{\EE XI(W\le q_\la(W))}%
{\EE XI(X\le q_\la(X))}.
$$
This is the same as the expression for $c_\la(X;W)$ with
$I(X\le q_\la(X))$ being replaced by~$X$.
However, an essential difference between~$\kappa$ and
the standard tail correlation coefficient is that
$\kappa$ is not symmetric in~$X,W$.

Let us also remark that, for the Weighted V@R,
$\kappa(X;W)$ remains unchanged under monotonic
transformations of~$W$, i.e.
$\kappa(X;W)=\kappa(X,f(W))$, where $f$ is a strictly
increasing function (this is seen from
Example~\ref{RC6}~(v)).

Let us study some basic properties of~$\kappa$.
Two propositions below correspond to two extremes:
no correlation and complete correlation.

\begin{Proposition}
\label{RC9}
Let $X,W\in L^1$ and suppose that $W$ has a
continuous distribution. Then
$u_\la^c(X;W)=0$ for any $\la\in(0,1]$
if and only if $\EE(X\cond W)=0$.
\end{Proposition}

\Proof
Set $Z_\la=\la^{-1}I(W\le q_\la(W))$.
According to Example~\ref{RC6}~(v),
$$
u_\la^c(X;W)
=\la^{-1}\EE XI(W\le q_\la(W))
=\la^{-1}\int_{(-\infty,q_\la(W)]}g(w)\QQ(d w),
\quad \la\in(0,1],
$$
where $g(w)=\EE(X\cond W=w)$ and $\QQ=\Law W$.
Now, the result is obvious.\End

\skm
Recall that random variables $X$ and $W$ are called
\textit{comonotone} if
$(X(\omega_2)-X(\omega_1))(W(\omega_2)-W(\omega_1))\ge0$
for $\PP\times\PP$-a.e. $\omega_1,\omega_2$.

\begin{Proposition}
\label{RC10}
Suppose that the support of~$\mu$ is $[0,1]$.
Let $X,W\in L^1(\DDD_\mu)$.
Then $u_\mu^c(X;W)=u_\mu(X)$ if and only if $X$ and $W$ are
comonotone.
\end{Proposition}

\Proof
Let us prove the ``only if'' part.
The map $\DDD_\mu\ni Z\mapsto\EE WZ$ is continuous with
respect to the weak topology $\si(L^1,L^\infty)$.
Hence, the set $\EX_{\DDD_\mu}(W)$ is weakly closed.
An application of the Hahn--Banach theorem shows that
$\EX_{\DDD_\mu}(W)$ is $L^1$-closed.
By Proposition~\ref{RC2}, there exists
$Z\in\EX_{\DDD_\mu}(W)$ such that $\EE XZ=u^c(X;W)$.
According to~\cite[Th.~4.4]{C05e}, there exists a
jointly measurable function
$(Z(\la,\omega);\la\in(0,1],\omega\in\Omega)$
such that $Z=\int_{(0,1]}Z_\la\mu(d\la)$ with
$Z_\la\in\DDD_\la$ for any~$\la$. Then
$$
\int_{(0,1]}\EE WZ_\la\mu(d\la)
=\EE WZ
=u_\mu(W)
=\int_{(0,1]}u_\la(W)\mu(d\la),
$$
and it follows that $Z_\la\in\EX_{\DDD_\la}(W)$ for
$\mu$-a.e.~$\la$. Furthermore,
$$
\int_{(0,1]}\EE XZ_\la\mu(d\la)
=\EE XZ
=u_\mu^c(X;W)
=u_\mu(X)
=\int_{(0,1]}u_\la(X)\mu(d\la),
$$
and it follows that $Z_\la\in\EX_{\DDD_\la}(X)$ for
$\mu$-a.e.~$\la$.
Thus, $\EX_{\DDD_\la}(W)\cap\EX_{\DDD_\la}(X)\ne\emp$
for $\mu$-a.e.~$\la$.
Using~\eqref{rc1}, we get
\begin{equation}
\label{rc7}
\PP(X>q_\la(X),W<q_\la(W))
=\PP(X<q_\la(X),W>q_\la(W))
=0
\end{equation}
for $\mu$-a.e.~$\la$. As the functions $q_\la(X)$ and
$q_\la(W)$ are right-continuous in~$\la$ and the support
of~$\mu$ is $[0,1]$, we deduce that~\eqref{rc7}
is satisfied for every~$\la\in(0,1]$.
From this it is easy to deduce that
$\PP((X,W)\in f((0,1]))=1$, where
$f(\la)=(q_\la(X),q_\la(W))$.
Thus, $X$ and $W$ are comonotone.

Let us prove the ``if'' part.
By~\cite[Lem.~4.83]{FS04}, there exists a random variable~$\xi$
and increasing functions $f,g$ such that
$X=f(\xi)$, $W=g(\xi)$. Set
$$
Z_\la=\la^{-1}I(\xi<q_\la(\xi))+cI(\xi=q_\la(\xi)),
$$
where $c$ is the constant such that $\EE Z_\la=1$.
According to~\cite[Th.~4.4]{C05e},
$Z:=\int_{(0,1]}Z_\la\mu(d\la)\in\DDD_\mu$.
It is clear that $Z_\la\in\EX_{\DDD_\la}(X)\cap\EX_{\DDD_\la}(W)$.
Hence,
\begin{align*}
&\EE WZ
=\int_{(0,1]}\EE WZ_\la\mu(d\la)
=\int_{(0,1]}u_\la(W)\mu(d\la)
=u_\mu(W),\\[1mm]
&\EE XZ
=\int_{(0,1]}\EE XZ_\la\mu(d\la)
=\int_{(0,1]}u_\la(X)\mu(d\la)
=u_\mu(X).
\end{align*}
As a result, $Z\in\EX_{\DDD_\mu}(W)$ and
$u^c(X;W)\le u_\mu(X)$. Since the reverse inequality is
obvious, we get $u_\mu^c(X;W)=u_\mu(X)$.\End

\skm
\textbf{5. Empirical estimation.}
In order to estimate empirically Alpha V@R contribution
of~$X$ to~$W$ with $\al\in\N$, one should first choose
the number of trials $K\in\N$ and generate independent
draws $(x_{kl},w_{kl};k=1,\dots,K,\,l=1,\dots,\al)$
of~$(X,W)$ using one of procedures described at the end
of Section~\ref{CRM}.
Having generated $x_{kl},w_{kl}$, one should calculate
the array
$$
l_k=\argmin_{l=1,\dots,\al}w_{kl},\quad k=1,\dots,K.
$$
According to Example~\ref{RC6}~(vii), an empirical
estimate of $\rho_\al^c(X;W)$ is provided by
$$
\rho_{\text{e}}^c(X;W)=-\frac{1}{K}\sum_{k=1}^K x_{kl_k}.
$$

In order to estimate Beta V@R contribution
with $\al,\be\in\N$, one should generate $x_{kl},w_{kl}$
similarly. Let $l_{k1},\dots,l_{k\be}$ be the numbers
$l\in\{1,\dots,\al\}$ such that the corresponding
$w_{kl}$ stand at the first $\be$ places
(in the increasing order) among
$w_{k1},\dots,w_{k\al}$.
According to Example~\ref{RC6}~(vi), an empirical
estimate of $\rho_{\al,\be}^c(X;W)$ is provided by
$$
\rho_{\text{e}}^c(X;W)
=-\frac{1}{K\be}\sum_{k=1}^K\sum_{i=1}^\be x_{kl_{ki}}.
$$

In order to estimate Weighted V@R contribution,
one should fix a data set
$(x_1,w_1),\dots,(x_T,w_T)$ and a measure~$\nu$ on
$\{1,\dots,T\}$ using one of the procedures described
at the end of Section~\ref{CRM}.
Let $w_{(1)},\dots,w_{(T)}$ be the values $w_1,\dots,w_T$
in the increasing order (we assume that all the $w_t$
are different).
Define $n(t)$ through the equality $w_{(t)}=w_{n(t)}$.
According to Example~\ref{RC6}~(iv), an empirical
estimate of $\rho_\mu^c(X;W)$ is provided by
$$
\rho_\text{e}^c(X;W)=
-\sum_{t=1}^T x_{n(t)}\int_{z_{t-1}}^{z_t}\psi_\mu(x)dx,
$$
where $z_t=\sum_{i=1}^T\nu\{n(i)\}$ and $\psi_\mu$
is given by~\eqref{crm7}.

%========================================================
\section{Factor Risk Contribution}
\label{FRC}

\textbf{1. Factor risk contribution.}
Let $(\Omega,\F,\PP)$ be a probability space
and $u$ be a coherent utility with the
determining set~$\DDD$.
Let $Y$ be a random variable (resp., random vector)
meaning the increment of some market factor
(resp., factors) over the unit time period and $W$ be
a random variable meaning the P\&L
produced by some portfolio over the unit time period.

\begin{Definition}\rm
\label{FRC1}
The \textit{factor utility contribution} is
$$
u^{fc}(X;Y;W)=\inf_{\QQ\in\EX_{\EE(\DDD\cond Y)}(W)}\EE_\QQ X,
\quad X\in L^0.
$$
The \textit{factor risk contribution} is
$\rho^{fc}(X;Y;W)=-u^{fc}(X;Y;W)$.
\end{Definition}

The function $u^{fc}(\,\cdot\,;Y;W)$ is a coherent utility
provided that $\EX_{\EE(\DDD\cond Y)}(W)\ne\emp$.

\begin{Proposition}
\label{FRC2}
If $\DDD$ is $L^1$-closed and uniformly integrable,
while $W\in L_s^1(\EE(\DDD\cond Y))$, then
$\EX_{\EE(\DDD\cond Y)}(W)\ne\emp$.
\end{Proposition}

\Proof
By Remark~(iii) following Definition~\ref{FR2},
$\EE(\DDD\cond Y)$ is $L^1$-closed and uniformly
integrable. Now, the result follows from
Proposition~\ref{RC2}.\End

\begin{Theorem}
\label{FRC3}
If $X,W\in L_s^1(\DDD)\cap L_s^1(\EE(\DDD\cond Y))\cap L^1$,
then
$$
u^{fc}(X;Y;W)=u^c(\EE(X\cond Y);\EE(W\cond Y)).
$$
\end{Theorem}

\Proof
By Lemma~\ref{FR4},
$$
\EE W\EE(Z\cond Y)=\EE\EE(W\cond Y)Z,\quad Z\in\DDD.
$$
Thus,
$$
\EE(Z\cond Y)\in\EX_{\EE(\DDD\cond Y)}(W)
\;\Lea\;Z\in\EX_\DDD(\EE(W\cond Y)),
$$
which means that
$$
\EX_{\EE(\DDD\cond Y)}(W)
=\EE\bigl(\EX_\DDD(\EE(W\cond Y))\cond Y\bigr).
$$
One more application of Lemma~\ref{FR4} yields
$$
\EE X\EE(Z\cond Y)=\EE\EE(X\cond Y)Z,\quad Z\in\DDD.
$$
As a result,
\begin{align*}
u^{fc}(X;Y;W)
&=\inf_{Z\in\EX_\DDD(\EE(W\cond Y))}\EE X\EE(Z\cond Y)\\
&=\inf_{Z\in\EX_\DDD(\EE(W\cond Y))}\EE\EE(X\cond Y)Z\\
&=u^c(\EE(X\cond Y);\EE(W\cond Y)).
\end{align*}

\vspace{-8mm}
\hfill$\Box$

\vspace{2mm}
\begin{Corollary}
\label{FRC4}
If $\DDD$ is $L^1$-closed and uniformly integrable,
while $X,W\in L_s^1(\DDD)\cap L_s^1(\EE(\DDD\cond Y))\cap L^1$,
then
$$
u^{fc}(X;Y;W)=\lim_{\eps\da0}\eps^{-1}
(u^f(W+\eps X;Y)-u^f(W;Y)).
$$
\end{Corollary}

\Proof
Applying successively Theorems~\ref{FRC3},
\ref{RC5}, and~\ref{FR3}, we get
\begin{align*}
u^{fc}(X;Y;W)
&=u^c(\EE(X\cond Y);\EE(W\cond Y))\\
&=\lim_{\eps\da0}\eps^{-1}(u(\EE(W\cond Y)
+\eps\EE(X\cond Y))-u(\EE(W\cond Y)))\\
&=\lim_{\eps\da0}\eps^{-1}(u^f(W+\eps X;Y)-u^f(W;Y))
\end{align*}
(in order to apply Theorem~\ref{RC5} in the second
equality, we need to check that
$\EE(X\cond Y)\in L_s^1(\DDD)$ and
$\EE(W\cond Y)\in L_s^1(\DDD)$; this is done by the same
argument as in the proof of Lemma~\ref{FR4}).\End

\skm
\Remark
If $(\Omega,\F,\PP)$ is atomless and $u$ is law invariant,
then, by Corollary~A.2, the integrability condition on
$X,W$ in the above statements can be replaced by a weaker one:
$X,W\in L_s^1(\DDD)$.\End

\begin{Example}\rm
\label{FRC5}
If $W$ is a constant,
then $\EX_{\EE(\DDD\cond Y)}(W)=\EE(\DDD\cond Y)$,
so that $u^{fc}(X;Y;W)=u^f(X;Y)$.

\textbf{(ii)} If $X=\la W$ with $\la\in\R_+$,
then $u^{fc}(X;Y;W)=\la u^f(W;Y)$
provided that $\EX_{\EE(\DDD\cond Y)}(W)\ne\emp$.

\textbf{(iii)} Let $u$ be law invariant and
$(X,Y,W)=(X,Y^1,\dots,Y^M,W)$ be a non-degenerate
Gaussian random vector such that each of its components
belongs to $L_s^1(\DDD)$. Let $C$ denote the covariance
matrix of~$Y$ and set $a^m=\cov(X,Y^m)$,
$b^m=\cov(W,Y^m)$, $\wl X=X-\EE X$, $\wl Y=Y-\EE Y$,
$\wl W=W-\EE W$.
We have (cf. Example~\ref{FR6})
$\EE(\wl X\cond\wl Y)=\lb C^{-1}a,\wl Y\rb$,
$\EE(\wl W\cond\wl Y)=\lb C^{-1}b,\wl Y\rb$,
so that, by Theorem~\ref{FRC3} and
Example~\ref{RC6}~(iii),
\begin{align*}
u^{fc}(X;Y;W)
&=\EE X+u^{fc}(\wl X;\wl Y;\wl W)\\
&=\EE X+u^c(\lb C^{-1}a,\wl Y\rb;
\lb C^{-1}b,\wl Y\rb)\\[1mm]
&=\EE X-\ga\,\frac{\cov(\lb C^{-1}a,\wl Y\rb,
\lb C^{-1}b,\wl Y\rb)}%
{(\var\lb C^{-1}b,\wl Y\rb)^{1/2}}\\[1mm]
&=\EE X-\ga\,\frac{\lb C^{-1}a,b\rb}%
{\lb C^{-1}b,b\rb^{1/2}}.
\end{align*}
In particular, if $Y$ is one-dimensional, then
$$
u^{fc}(X;Y;W)=\EE X-\ga\,\frac{\cov(X,Y)}%
{(\var Y)^{1/2}}\,\sgn\cov(W,Y).
$$
If moreover, $\cov(X,Y)>0$ and $\cov(W,Y)>0$, then, recalling
Example~\ref{FR6} and Example~\ref{RC6}~(iii), we get
$$
u^f(X;Y)
=u^c(X;Y)
=u^{fc}(X;Y;W)
=\EE X-\ga\,\frac{\cov(X,Y)}{(\var Y)^{1/2}}.
$$

\vspace{-8mm}
\hfill$\Box$
\end{Example}

\vspace{5mm}
\textbf{2. Empirical estimation.}
In view of Theorem~\ref{FRC3}, the empirical estimation
of $\rho^{fc}(X;Y;W)$ reduces to finding the functions
$f(y)=\EE(X\cond Y=y)$, $g(y)=\EE(W\cond Y=y)$
and then applying the procedures described at the end of
Section~\ref{RC} to $x_{kl}:=f(y_{kl})$,
$w_{kl}=g(y_{kl})$.

If $Y$ is one-dimensional, one can create the data for~$Y$
using the time change or the scaling procedures described
at the end of Section~\ref{FR}.

%======================================================
\section{Optimal Risk Sharing}
\label{ORS}

\textbf{1. Problem.}
Let $(\Omega,\F,\PP)$ be a probability space and
$u^1,\dots,u^M$ be coherent utilities with the
determining sets $\DDD^1,\dots,\DDD^M$.
We assume that each $\DDD^m$ is $L^1$-closed and
uniformly integrable.
Suppose there is a firm consisting of~$N$ desks,
and the $n$-th desk can invest into the assets
that produce P\&Ls $(X^{n1},\dots,X^{nd^n})$.
We assume that $X^{nk}\in L_s^1(\DDD^m)$ for any $k,n,m$.
The set of P\&Ls that the $n$-th desk can produce over a
unit time period is $\{\lb h^n X^n\rb:h^n\in H^n\}$, where
$X^n=(X^{n1},\dots,X^{nd^n})$ and $H^n$ is a convex
subset of~$\R^{d^n}$ with a non-empty interior meaning
the constraint on the portfolio of the $n$-th desk.
We assume that $u^m(\lb h^n,X^n\rb)<0$ for any
$h^n\in H^n\setminus\{0\}$, which means that any possible
trade has a strictly positive risk.
Let $E^n\in\R^{d^n}$ be the vector of rewards for the
assets available to the $n$-th desk.
This is the vector of subjective assessments by the
$n$-th desk of the profitability of the assets
$X^{n1},\dots,X^{nd^n}$.

We will consider the following optimization problem
for the whole firm:
\begin{equation}
\label{ors1}
\begin{cases}
\sum_n\lb h^n,E^n\rb\longrightarrow\max,\\
h^n\in H^n,\;n=1,\dots,N,\\[0.5mm]
\rho^m\bigl(\sum_n\lb h^n,X^n\rb\bigr)\le c^m,\;m=1,\dots,M,
\end{cases}
\end{equation}
where $c^1,\dots,c^M\in(0,\infty)$ are fixed risk limits.
This problem is a generalization of~\eqref{po1}, which
might be considered as the optimization
problem for a separate desk.
We will not try to solve~\eqref{ors1} for the following
reason: if this problem admitted a solution that can be
implemented in practice, this would mean that the central
management is able to optimize the firm's portfolio
and there would be no need for the existence of
separate desks.
Instead of trying to solve~\eqref{ors1}, we will study
the following question:
\textit{Is it possible to decentralize~\eqref{ors1},
i.e. to create for the desks conditions such that
the global optimum is achieved
when each desk acts optimally?}

\skm
\textbf{Hypothesis~1.}
There exist $c^{nm}$ such that if $h_*^n$
satisfy the conditions
\begin{mitemize}
\item[1.] we have
$$
\rho^m\Bigl(\sum_{n=1}^N \lb h_*^n,X^n\rb\Bigr)\le c^m,
\quad m=1,\dots,M,
$$
and the equality is attained at least for one~$m$;
\item[2.] for each~$n$, the vector $h_*^n$ solves the problem
\begin{equation}
\label{ors2}
\begin{cases}
\lb h^n,E^n\rb\longrightarrow\max,\\
h^n\in H^n,\\
\rho^m(\lb h^n,X^n\rb)\le c^{nm},\;m=1,\dots,M,
\end{cases}
\end{equation}
\end{mitemize}
then $(h_*^1,\dots,h_*^N)$ solves~\eqref{ors1}.

\skm
This hypothesis is wrong as shown by the example below.

\begin{Example}\rm
\label{ORS1}
Let $M=1$, $N=2$, $u$ be a law invariant coherent utility
that is finite on Gaussian random variables,
$H^1=\R$, $H^2=\R^2$, and $(X^1,X^{21},X^{22})$
have a jointly Gaussian distribution
with a non-degenerate covariance matrix~$C$.
It was shown in~\cite[Subsect.~2.2]{C062} that the
solution of~\eqref{ors1} has the form
$$
(h_*^1,h_*^{21},h_*^{22})
=const\times C^{-1}(E^1,E^{21},E^{22}).
$$
Furthermore, the solution of~\eqref{ors2} with $n=2$
has the form
$$
(\wt h_*^{21},\wt h_*^{22})=const\times
\wt C^{-1}(E^{21},E^{22}),
$$
where $\wt C$ is the covariance matrix of
$(X^{21},X^{22})$ (the constant here depends on $c^{nm}$).
It is clear that there need not exist $c^{nm}$ such
that $(\wt h_*^{21},\wt h_*^{22})=(h_*^{21},h_*^{22})$.\End
\end{Example}

\textbf{2. Limits on risk contribution.}
The reason why Hypothesis~1 is wrong
is that in general
$$
\rho\Bigl(\sum_{n=1}^N\lb h^n,X^n\rb\Bigr)
<\sum_{n=1}^N\rho(\lb h^n,X^n\rb).
$$
On the other hand, we typically have
$$
\rho\Bigl(\sum_{n=1}^N\lb h^n,X^n\rb\Bigr)
=\sum_{n=1}^N\rho^c\Bigl(\lb h^n,X^n\rb;
\sum_{n=1}^N\lb h^n,X^n\rb\Bigr)
$$
(see~\eqref{rc6}). This gives rise to

\skm
\textbf{Hypothesis~2.}
Let $c^{nm}\in\R_+$ be such that
$\sum_n c^{nm}=c^m$ for each~$m$.
If $h_*^n$ satisfy the conditions
\begin{mitemize}
\item[1.] we have
$$
\rho^m\Bigl(\sum_{n=1}^N \lb h_*^n,X^n\rb\Bigr)\le c^m,
\quad m=1,\dots,M,
$$
and the equality is attained at least for one~$m$;
\item[2.] for each~$n$, the vector $h_*^n$ solves the problem
$$
\begin{cases}
\lb h^n,E^n\rb\longrightarrow\max,\\
h^n\in H^n,\\
(\rho^m)^c\bigl(\lb h^n,X^n\rb;\sum_n\lb h_*^n,X^n\rb\bigr)
\le c^{nm},\;m=1,\dots,M,
\end{cases}
$$
\end{mitemize}
then $(h_*^1,\dots,h_*^N)$ solves~\eqref{ors1}.

\skm
This hypothesis is also wrong as shown by the example
below.

\begin{Example}\rm
\label{ORS2}
Let $M=1$, $N=2$, $u$ be a law invariant coherent
utility, $H^1=H^2=\R$, $X^1,X^2\in L_s^1(\DDD)$
be jointly Gaussian with a non-degenerate covariance
matrix, and $E^1=E^2=1$.
Take arbitrary $h_*^1,h_*^2\in(0,\infty)$ such that
$\rho(h_*^1X^1+h_*^2X^2)=c$ and
$$
\EE X^n<\ga\,\frac{\cov(X^n,h_*^1X^1+h_*^2X^2)}%
{(\var(h_*^1X^1+h_*^2 X^2))^{1/2}},\quad n=1,2,
$$
where $\ga$ is provided by Example~\ref{CRM11}~(i). Set
$$
c^n:=\rho^c(h_*^n X^n;h_*^1X^1+h_*^2X^2)
=-h_*^n\EE X^n+h_*^n\ga\,\frac{\cov(X^1,h_*^1X^1+h_*^2X^2)}%
{(\var(h_*^1X^1+h_*^2 X^2))^{1/2}},\quad n=1,2
$$
(we used Example~\ref{RC6}~(iii)).
Obviously, each $h_*^n$ solves~\eqref{ors2}.
On the other hand, $(h_*^1,h_*^2)$ need not be optimal
for~\eqref{ors1}.\End
\end{Example}

\textbf{3. Risk trading.}
The reason why Hypothesis~2 is wrong is that if one unit
is more profitable than another, but it obtains a
lower risk limit, then the global optimum cannot be
achieved. In fact, if we replace fixed~$c^{nm}$ in
Hypothesis~2 by the assumption ``there exists $c^{nm}$...'',
then (as follows from Theorem~\ref{ORS3}) the hypothesis
becomes true. But this leaves open the problem of
finding~$c^{nm}$. Instead of trying to solve this
problem, we will take another path. Let us assume that
the desks are allowed to trade their risk limits, i.e.
they establish themselves (through the supply--demand
equilibrium) the price~$\al^m$ for the $m$-th risk limit,
so that if the $n$-th desk buys from the $n'$-th desk
$a$ units of the $m$-th risk limit, then the $n$-th
desk pays the $n'$-th desk the amount $a\al^m$,
the $n$-th desk raises its $m$-th risk limit by
the amount~$a$,
and the $n'$-th desk lowers its $m$-th risk limit by
the amount~$a$.

\skm
\textbf{Hypothesis~3.}
Let $c^{nm}\in\R_+$ be such that $\sum_n c^{nm}=c^m$ for
each~$m$. If $h_*^n\in H^n$, $a_*^{nm}\in\R$, and
$\al_*^m\in\R_+$ satisfy the conditions
\begin{mitemize}
\item[1.] we have
$$
\sum_{n=1}^N a_*^{nm}=0,\quad m=1,\dots,M;
$$
\item[2.] we have
\begin{equation}
\label{ors3}
\rho^m\Bigl(\sum_{n=1}^N\lb h_*^n,X^n\rb\Bigr)\le c^m,
\quad m=1,\dots,M,
\end{equation}
and the equality is attained at least for one~$m$;
\item[3.] $\al_*^m=0$ for all~$m$ such that
inequality~\eqref{ors3} is strict;
\item[4.] for each~$n$, the vectors $h_*^n$ and
$a_*^n=(a_*^{nm};m=1,\dots,M)$ solve the problem
$$
\begin{cases}
\lb h^n,E^n\rb-\lb a^n,\al_*\rb\longrightarrow\max,\\
h^n\in H^n,\;a^n\in\R^m,\\
(\rho^m)^c\bigl(\lb h^n,X^n\rb;\sum_n\lb h_*^n,X^n\rb\bigr)
\le c^{nm}+a^{nm},\;m=1,\dots,M,
\end{cases}
$$
\end{mitemize}
then $(h_*^1,\dots,h_*^N)$ solves~\eqref{ors1}.

\skm
This hypothesis is true (under minor technical
assumptions) as shown by the theorem below.

\begin{Theorem}
\label{ORS3}
Let $c^{nm}\in\R_+$ be such that $\sum_n c^{nm}=c^m$ for
each~$m$.
Let $h_*^n$ belong to the interior of~$H^n$ and assume
that, for each~$m$, the set
$\EX_{\DDD^m}\bigl(\sum_n\lb h_*^n,X^n\rb\bigr)$
consists of a unique measure~$\QQ^m$.
Then the following conditions are equivalent:
\begin{mitemize}
\item[\rm(i)] $h_*^1,\dots,h_*^N$ solve~\eqref{ors1}{\rm;}
\item[\rm(ii)] there exist $a_*^{nm}\in\R$ and
$\al_*^m\in\R_+$ that satisfy the conditions of
Hypothesis~3;
\item[\rm(iii)] there exist $\al_*^m\in\R_+$ satisfying
conditions 2, 3 of Hypothesis~3 and such that
$$
E^n=-\sum_{m=1}^M\al_*^m\EE_{\QQ^m}X^n,\quad n=1,\dots,N.
$$
\end{mitemize}
Moreover, the sets of possible $\al_*$ in {\rm(ii)} and
{\rm(iii)} are the same.
\end{Theorem}

\Remark
The theorem shows that the system finds the optimum
regardless of the allocation $c^{nm}$ of risk limits.
(The resulting rewards
$\lb h_*^n,E^n\rb-\lb a_*^n,\al_*\rb$ depend
on~$c^{nm}$, but their sum does not depend on~$c^{nm}$.)
It is also clear from~(iii) that the equilibrium prices
$\al_*^m$ of risk limits do not depend on $c^{nm}$.\End

\skm
\texttt{Proof of Theorem~\ref{ORS3}.}
(i)$\Rightarrow$(iii) Set
\begin{align*}
J&=\Bigl\{m:\rho^m\Bigl(\sum_{n=1}^N\lb h_*^n,
X^n\rb\Bigr)=c^m\Bigr\},\\
K&=\Bigl\{-\sum_{m\in J}\al^m\EE_{\QQ^m}(X^1,\dots,X^N):
\al^m\in\R_+\Bigr\},
\end{align*}
so that $K$ is a convex closed cone in~$\R^d$, where
$d=d^1+\dots+d^N$.
Suppose that $E:=(E^1,\dots,E^N)\notin K$. By the Hahn--Banach theorem,
there exists $h=(h^1,\dots,h^N)\in\R^d$ such that
$\sup_{x\in K}\lb h,x\rb\le0<\lb h,E\rb$.
Consider $h^n(\eps)=h_*^n+\eps h^n$.
There exists $\de>0$ such that, for $\eps\in(0,\de)$,
we have: $h^n(\eps)\in H^n$ for any~$n$ and
$\rho^m\bigl(\sum_n\lb h^n,X^n\rb\bigr)<c^m$ for any
$m\notin J$. Set
\begin{align*}
f(\eps)&=\sum_{n=1}^N\lb h^n(\eps),E^n\rb,\\
g^m(\eps)&=\rho^m\Bigl(\sum_{n=1}^N\lb h^n(\eps),X^n\rb\Bigr),
\quad m\in J.
\end{align*}
Then
$$
\frac{d}{d\eps}\Bigl|_{\eps=0}f(\eps)
=\sum_{n=1}^N\lb h^n,E^n\rb
=\lb h,E\rb
>0,
$$
and, by Theorem~\ref{RC5},
$$
\frac{d}{d\eps}\Bigl|_{\eps=0} g^m(\eps)
=-\EE_{\QQ^m}\sum_{n=1}^N\lb h^n,X^n\rb
\le\sup_{x\in K}\lb h,x\rb
\le0,
\quad m\in J.
$$
Due to our assumptions, $g^m(0)>0$.
Then we can find $\eps\in(0,\de)$ and $\be\in(0,1)$
such that $\be h^n(\eps)\in H^n$ for each~$n$,
$\be f(\eps)>f(0)$, and $\be g^m(\eps)\le g^m(0)$ for any
$m\in J$. This means that
$$
\sum_{n=1}^N\lb\be h^n(\eps),E^n\rb
>\sum_{n=1}^N\lb h_*^n,E^n\rb
$$
and
$$
\rho^m\Bigl(\sum_{n=1}^N\lb \be h^n(\eps),X^n\rb\Bigr)\le c^m,
\quad m=1,\dots,M,
$$
which is a contradiction.
As a result, $E\in K$, which is the desired statement.

(iii)$\Rightarrow$(ii)
It follows from the inequality
$$
-\sum_{n=1}^N\EE_{\QQ^m}\lb h_*^n, X^n\rb
=\rho^m\Bigl(\sum_{n=1}^N\lb h_*^n,X^n\rb\Bigr)
\le c^m
=\sum_{n=1}^N c^{nm},\quad m=1,\dots,M
$$
that we can choose $a_*^{nm}$ in such a way that
$\sum_n a_*^{nm}=0$ for any~$m$ and
$a_*^{nm}\ge-\EE_{\QQ^m}\lb h_*^n,X^n\rb-c^{nm}$
for any $n,m$. Then
$$
(\rho^m)^c\Bigl(\lb h_*^n,X^n\rb;
\sum_{n=1}^N\lb h_*^n,X^n\rb\Bigr)
=-\EE_{\QQ^m}\lb h_*^n,X^n\rb
\le c^{nm}+a_*^{nm},\quad m=1,\dots,M.
$$
Clearly, for $m\in J$, this inequality is the equality.
Assume that there exist~$n$ and $h^n\in H^n$,
$a^n\in\R^m$ such that
\begin{gather*}
\lb h^n,E^n\rb-\lb a^n,\al_*\rb
>\lb h_*^n,E^n\rb-\lb a_*^n,\al_*\rb,\\
(\rho^m)^c\Bigl(\lb h^n,X^n\rb;
\sum_{n=1}^N\lb h_*^n,X^n\rb\Bigr)\le c^{nm}+a^{nm},
\quad m=1,\dots,M.
\end{gather*}
This means that, for $\Delta h^n=h^n-h_*^n$ and
$\Delta a^n=a^n-a_*^n$, we have
\begin{gather*}
\lb \Delta h^n,E^n\rb>\lb\Delta a^n,\al_*\rb,\\
\EE_{\QQ^m}\lb\Delta h^n,X^n\rb\ge-\Delta a^{nm},
\quad m\in J.
\end{gather*}
This leads to the inequality
$$
\lb\Delta h^n,E^n\rb+\sum_{m\in J}\al_*^m
\EE_{\QQ^m}\lb\Delta h^n,X^n\rb>0,
$$
which is a contradiction.

(ii)$\Rightarrow$(i)
Suppose that there exist $h^n\in H^n$ such that
\begin{gather*}
\sum_{n=1}^N\lb h^n,E^n\rb
>\sum_{n=1}^N\lb h_*^n,E^n\rb,\\
\rho^m\Bigl(\sum_{n=1}^N\lb h^n,X^n\rb\Bigr)\le c^m,
\quad m=1,\dots,M.
\end{gather*}
Then $h^n(\eps):=h_*^n+\eps(h^n-h_*^n)\in H^n$
for any $\eps\in(0,1)$.
For $\De h^n=h^n-h_*^n$, we have
$\sum_n\lb\De h^n,E^n\rb>0$, and, due to the convexity
of the function
$\eps\mapsto\rho^m\bigl(\sum_n\lb h^n(\eps),X^n\rb\bigr)$,
$$
\sum_{n=1}^N\EE_{\QQ^m}\lb\De h^n,X^n\rb
=-\frac{d}{d\eps}\Bigl|_{\eps=0}\rho^m\Bigl(
\sum_{n=1}^N\lb h^n(\eps),X^n\rb\Bigr)
\ge0,\quad m\in J.
$$
This means that there exists~$n$ such that
$$
\lb\De h^n,E^n\rb
>-\sum_{m\in J}\al_*^m\EE_{\QQ^m}\lb\De h^n,X^n\rb.
$$
For this~$n$, set
$a^{nm}=a_*^{nm}-\EE_{\QQ^m}\lb\De h^n,X^n\rb$.
Then
$$
\lb h^n,E^n\rb-\lb a^n,\al_*\rb
>\lb h_*^n,E^n\rb-\lb a_*^n,\al_*\rb.
$$
Furthermore,
\begin{align*}
(\rho^m)^c\Bigl(\lb h^n,X^n\rb;\sum_{n=1}^N
\lb h_*^n,X^n\rb\Bigr)
&=-\EE_{\QQ^m}\lb h_*^n,X^n\rb
-\EE_{\QQ^m}\lb\Delta h^n,X^n\rb\\[-3mm]
&=(\rho^m)^c\Bigl(\lb h_*^n,X^n\rb;\sum_{n=1}^N
\lb h_*^n,X^n\rb\Bigr)+a^{nm}-a_*^{nm}\\
&\le c^{nm}+a^{nm},\quad m=1,\dots,M,
\end{align*}
which is a contradiction.

In order to complete the proof, it is sufficient to
show that any $\al_*$ satisfying~(ii) also satisfies~(iii).
This is obvious.\End

%======================================================
\section{Summary and Conclusion}
\label{SC}

\textbf{1. Risk measurement.}
Consider a firm whose portfolio has the
form $W=\sum_{n=1}^N W^n$.
Here $W$ is the P\&L produced by the portfolio
over the unit time period~$\Delta$, which is used as
the basis for risk measurement (typically it is one day);
$W^n$ is the P\&L of the $n$-th asset.
Let $X$ be the P\&L produced by some additional asset
or portfolio over the same period.

The empirical procedure to assess Alpha V@R risk
of~$W$ and the risk contribution of~$X$ to~$W$ is:
\begin{mitemize}
\item[\bf 1.] Fix $\al\in\N$.
\item[\bf 2.] Choose a probability measure~$\nu$ on
the set of natural numbers,
choose the number of trials $K\in\N$,
and generate independent draws
$(t_{kl};k=1,\dots,K,\,l=1,\dots,\al)$ from the
distribution~$\nu$. A natural choice for~$\nu$ is a
geometric distribution.
\item[\bf 3.] Calculate the array
$$
l_k=\argmin_{l=1,\dots,\al}\sum_{n=1}^N w_{kl}^n,
\quad k=1,\dots,K.
$$
Here $w_{kl}^n$ is the increment of the value of the
$n$-th asset in the portfolio over the time period
$[-t_{kl}\De,-(t_{kl}-1)\De]$ (the current time
instant is~0)\footnote{Of course, if the price of the $n$-th
asset is positive by its nature (for example, the
$n$-th asset is a share), then one can use a finer
procedure to determine $w_{kl}^n$,
i.e. $w_{kl}^n$ is the current value of the asset times
its relative increment over the period
$[-t_{kl}\De,-(t_{kl}-1)\De]$.}.
\item[\bf 4.] Calculate the empirical estimates of
$\rho_\al(W)$ and $\rho_\al^c(X;W)$ by the formulas:
\begin{align*}
&\rho_{\text{e}}(W)
=-\frac{1}{K}\sum_{k=1}^K\sum_{n=1}^N w_{kl_k}^n,\\
&\rho_{\text{e}}^c(X;W)
=-\frac{1}{K}\sum_{k=1}^K x_{kl_k}.
\end{align*}
Here $x_{kl}$ is the realization of~$X$
(i.e. the increment of the value of the corresponding
portfolio) over the time period
$[-t_{kl}\De,-(t_{kl}-1)\De]$.
\end{mitemize}

Note that the same arrays $(t_{kl})$ and $(l_k)$ are
used for different~$X$. If we measure risk on the daily
basis, steps~1--3 can be performed only once a day
(for example, in the night).
Thus, in the morning the central desk announces the arrays
$(t_{kl})$ and $(l_k)$, and when assessing the risk contribution of any
trade~$X$, any desk should simply take the realizations
of~$X$ over the corresponding intervals and insert
them into the formula for~$\rho_{\text{e}}^c$.
The number of operations required to generate
$(t_{kl})$ and $(l_k)$ is of order $\al N K$;
the number of operations required to calculate
$\rho_{\text{e}}^c$ is of order~$K$.
In particular, we need not order the data set
(which is needed for estimating V@R).

The above estimation procedure is completely non-linear.
Moreover, it is completely empirical as
it uses no model assumptions on the structure
of~$X$ and~$W$.

If $X=\sum_{j=1}^{J}X^j$ is itself a big portfolio
(for example, the P\&L produced by a desk of the firm),
then both the theoretical risk contribution and its
empirical estimate satisfy the linearity property:
$$
\rho^c\Bigl(\sum_{j=1}^J X^j;W\Bigr)
=\sum_{j=1}^J\rho^c(X^j;W).
$$

The above procedure can be combined with the bootstrap
technique, with obvious changes.

A similar procedure can be performed for Beta V@R.
The difference is that one should additionally
choose $\be\in\{1,\dots,\al-1\}$ and find the numbers
$l_{k1},\dots,l_{k\be}$ such that the corresponding
$\sum_n\!w_{kl}^n$ stand at the first $\be$ places
(in the increasing order) among
$\sum_n\!w_{k1}^n,\dots,\sum_n\!w_{k\al}^n$.
Then the empirical estimates of $\rho_{\al,\be}(W)$ and
$\rho_{\al,\be}^c(X;W)$ are provided by
\begin{align*}
&\rho_{\text{e}}(W)
=-\frac{1}{K\be}\sum_{k=1}^K\sum_{i=1}^\be\sum_{n=1}^N
w^n_{kl_{ki}},\\
&\rho_{\text{e}}^c(X;W)
=-\frac{1}{K\be}\sum_{k=1}^K\sum_{i=1}^\be x_{kl_{ki}}.
\end{align*}

\skb
\textbf{2. Factor risk measurement.}
The procedure to assess Alpha V@R factor risks of~$W$
and the factor risk contributions of~$X$ to~$W$ is:
\begin{mitemize}
\item[\bf 1.] Fix $\al\in\N$.
\item[\bf 2.] Choose the main market factors
$Y^1,\dots,Y^M$ affecting the risk of the portfolio.
Here $Y^m$ means the increment of the $m$-th factor
over the unit time period.
\item[\bf 3.] Create procedures for calculating the
functions $f^m(y)=\EE(X\cond Y^m=y)$ and
$g^{nm}(y)=\EE(W^n\cond Y^m=y)$.
\item[\bf 4.] Choose a probability measure~$\nu$ on the
set of natural numbers, choose the number of trials~$K$,
and generate independent draws
$(t_{kl};k=1,\dots,K,\,l=1,\dots,\al)$ from the
distribution~$\nu$.
\item[\bf 5.] Calculate the array
$$
l_k^m=\argmin_{l=1,\dots,\al}\sum_{n=1}^N
g^{nm}(y_{kl}^m),\quad k=1,\dots,K,\;m=1,\dots,M.
$$
Here $y_{kl}^m$ is the increment of the $m$-th factor
over the time period
$[-t_{kl}(\si^m)^2\De,-(t_{kl}-1)(\si^m)^2\De]$,
where $\si^m$ is the current volatility of the $m$-th
factor (for example, the implied volatility).
Instead of this time change procedure, one can use the
scaling procedure.
\item[\bf 6.] Calculate the empirical estimates of
$\rho_\al^f(W;Y^m)$ and $\rho_\al^{fc}(X;Y^m;W)$ by the formulas:
\begin{align*}
&\rho_{\text{e}}^f(W;Y^m)
=-\frac{1}{K}\sum_{k=1}^K\sum_{n=1}^N
g^{nm}\bigl(y_{kl_k^m}^m\bigr),\quad m=1,\dots,M,\\
&\rho_{\text{e}}^{fc}(X;Y^m;W)
=-\frac{1}{K}\sum_{k=1}^K f^m\bigl(y_{kl_k^m}^m\bigr),
\quad m=1,\dots,M.
\end{align*}
\end{mitemize}

All the pleasant features of risk estimates described above
remain true for factor risks.
Moreover, an important advantage of factor risks is
that we can take an arbitrarily large data set
$y_1,\dots,y_T$, while the joint data set
$w_t^n$ required for ordinary risk does not exist for
large portfolios.

A similar procedure can be performed for Beta V@R.
The difference is that one should additionally
choose $\be\in\{1,\dots,\al-1\}$ and find the numbers
$l_{k1}^m,\dots,l_{k\be}^m$ such that the corresponding
$\sum_n\!g^{nm}(y_{kl}^m)$ stand at the first $\be$ places
(in the increasing order) among
$\sum_n\!g^{nm}(y_{k1}^m),\dots,
\sum_n\!g^{nm}(y_{k\al}^m)$.
Then the empirical estimates of
$\rho_{\al,\be}^f(W;Y^m)$ and
$\rho_{\al,\be}^{fc}(X;Y^m;W)$ are provided by
\begin{align*}
&\rho_{\text{e}}^f(W;Y^m)
=-\frac{1}{K\be}\sum_{k=1}^K\sum_{i=1}^\be\sum_{n=1}^N
g^{nm}\bigl(y_{kl_{ki}^m}^m\bigr),\quad m=1,\dots,M,\\
&\rho_{\text{e}}^{fc}(X;Y^m;W)
=-\frac{1}{K\be}\sum_{k=1}^K\sum_{i=1}^\be
f^m\bigl(y_{kl_{ki}^m}^m\bigr),
\quad m=1,\dots,M.
\end{align*}

The functions $g^{nm}$ should not be recalculated every
day (if we measure risk on a daily basis);
once computed, these functions
might be used for rather a long period.
Of course, the collection of assets in the portfolio
changes each day, but the main part of this collection
remains the same, so that when the $(N+1)$-th asset
appears in the portfolio, we should calculate
$g^{N+1,m}$, $m=1,\dots,M$, but we should not recalculate
$g^{nm}$, $n=1,\dots,N$, $m=1,\dots,M$.

The factor risk estimation admits a multi-factor version:
instead of considering $M$ different factor risks,
we consider one risk driven by the multidimensional
factor $Y=(Y^1,\dots,Y^M)$.
Alternatively, we could split the factors into $M$ groups
so that the increment of the $m$-th group is a random
vector, which we still denote by~$Y^m$,
and split the portfolio into $M$ groups so that the risk
of the $m$-th group is driven mainly by the $m$-th group
of factors. Then when dealing with the $m$-th factor
risk we are considering only the $m$-th part of the
portfolio.
All the formulas remain the same with obvious changes.

\skm
\textbf{3. Comparison of various techniques.}
In Table~1, we compare different risk measurement
techniques proposed in this paper and several classical
risk measurement techniques (see~\cite[Sect.~6]{M02}
for their description).
Let us briefly describe this table.

All the techniques, except for one-factor risk
measurement, assess the overall risk.
As shown by Example~\ref{FR11}, the sum of one-factor
risks need not exceed the multi-factor risk, so that
the measurement of one-factor risks might be insufficient.

For the empirical risk estimation procedure described
above, one requires the joint time series $w_t^n$ for
all the assets in the portfolio.
However, different assets have different
durations, so that a joint time series might exist
only for small portfolios or portfolios consisting
of long-living assets only.
In contrast, all the other methods express the
values of the assets through the main factors, and for
the factors one can get arbitrarily large time series.

For one-factor risks
we can use the time change technique described above,
which enables one to react immediately to the volatility
changes. For parametric V@R and Monte Carlo V@R, this
can partially be done, but these methods require the
covariance matrix for different assets, for which we
should use historic data.

The speed of computations for one-factor coherent risks,
multi-factor coherent risks, and historic V@R is
approximately the same because in all these methods
the main part of computations consists in finding the
values $g^{nm}(y_{kl}^m)$.

The methods proposed in this paper deal with coherent risk.
The arguments of Section~\ref{CRM} show that in many
respects this is much wiser than the use of V@R.

All the methods of this paper admit simple
calculation of risk contributions.
For parametric V@R, this is also possible because
there exist explicit formulas for Gaussian V@R
contributions.
For Monte Carlo V@R, the possibility to calculate
risk contributions depends on the choice of the
probabilistic model.

All the methods, except for parametric V@R, are completely
non-linear.

Monte Carlo V@R might be performed both for Gaussian
and non-Gaussian models, but for the latter ones this
leads to a further decrease in the speed of computations.

The empirical risk estimation technique uses absolutely
no model assumptions. As for the factor risks and
historic V@R, there are no assumptions on the
probabilistic structure of~$Y^m$, but model assumptions
are used when calculating the functions $f^m$, $g^{nm}$.

As the conclusion, the best methods are:
one-factor and multi-factor coherent risk measurement.
It might be reasonable to use both of them simultaneously.

\vspace{5mm}
\begin{center}
%\begin{center}
\begin{tabular}{|p{33mm}||p{9mm}|p{9mm}|p{9mm}|%
p{9mm}|p{9mm}|p{9mm}|p{9mm}|p{9mm}|p{9mm}|}
\hline
&\rotatebox{90}{\vbox{
\hbox{Assessment of}
\hbox{the overall risk }}}&
\rotatebox{90}{\vbox{
\hbox{Availability}
\hbox{of data}}}&
\rotatebox{90}{\vbox{
\hbox{Ability to capture}
\hbox{volatility predictions }}}&
\rotatebox{90}{\vbox{
\hbox{Speed of}
\hbox{computations}}}&
\rotatebox{90}{\vbox{
\hbox{Coherence}
\hbox{properties}}}&
\rotatebox{90}{\vbox{
\hbox{Ability to calculate}
\hbox{risk contributions}}}&
\rotatebox{90}{\vbox{
\hbox{Ability to capture}
\hbox{non-linearity}}}&
\rotatebox{90}{\vbox{
\hbox{Ability to capture}
\hbox{non-normality}}}&
\rotatebox{90}{\vbox{
\hbox{Model}
\hbox{independence}}}\\
\hhline{|=||=|=|=|=\=|=|=|=|=|=|}
\rule{0mm}{5.5mm}Ordinary risk&
\hs\raisebox{-1.2mm}{\includegraphics{factor.7}}&
\hs\raisebox{-1.2mm}{\includegraphics{factor.8}}&
\hs\raisebox{-1.2mm}{\includegraphics{factor.9}}&
\hs\raisebox{-1.2mm}{\includegraphics{factor.7}}&
\hs\raisebox{-1.2mm}{\includegraphics{factor.7}}&
\hs\raisebox{-1.2mm}{\includegraphics{factor.7}}&
\hs\raisebox{-1.2mm}{\includegraphics{factor.7}}&
\hs\raisebox{-1.2mm}{\includegraphics{factor.7}}&
\hs\raisebox{-1.2mm}{\includegraphics{factor.7}}\\[2mm]
\hline
\rule{0mm}{5.5mm}Multi-factor risk&
\hs\raisebox{-1.2mm}{\includegraphics{factor.7}}&
\hs\raisebox{-1.2mm}{\includegraphics{factor.7}}&
\hs\raisebox{-1.2mm}{\includegraphics{factor.9}}&
\hs\raisebox{-1.2mm}{\includegraphics{factor.8}}&
\hs\raisebox{-1.2mm}{\includegraphics{factor.7}}&
\hs\raisebox{-1.2mm}{\includegraphics{factor.7}}&
\hs\raisebox{-1.2mm}{\includegraphics{factor.7}}&
\hs\raisebox{-1.2mm}{\includegraphics{factor.7}}&
\hs\raisebox{-1.2mm}{\includegraphics{factor.8}}\\[2mm]
\hline
\rule{0mm}{5.5mm}One-factor risk&
\hs\raisebox{-1.2mm}{\includegraphics{factor.9}}&
\hs\raisebox{-1.2mm}{\includegraphics{factor.7}}&
\hs\raisebox{-1.2mm}{\includegraphics{factor.7}}&
\hs\raisebox{-1.2mm}{\includegraphics{factor.8}}&
\hs\raisebox{-1.2mm}{\includegraphics{factor.7}}&
\hs\raisebox{-1.2mm}{\includegraphics{factor.7}}&
\hs\raisebox{-1.2mm}{\includegraphics{factor.7}}&
\hs\raisebox{-1.2mm}{\includegraphics{factor.7}}&
\hs\raisebox{-1.2mm}{\includegraphics{factor.8}}\\[2mm]
\hline
\rule{0mm}{5.5mm}Parametric V@R&
\hs\raisebox{-1.2mm}{\includegraphics{factor.7}}&
\hs\raisebox{-1.2mm}{\includegraphics{factor.7}}&
\hs\raisebox{-1.2mm}{\includegraphics{factor.8}}&
\hs\raisebox{-1.2mm}{\includegraphics{factor.7}}&
\hs\raisebox{-1.2mm}{\includegraphics{factor.9}}&
\hs\raisebox{-1.2mm}{\includegraphics{factor.7}}&
\hs\raisebox{-1.2mm}{\includegraphics{factor.9}}&
\hs\raisebox{-1.2mm}{\includegraphics{factor.9}}&
\hs\raisebox{-1.2mm}{\includegraphics{factor.9}}\\[2mm]
\hline
\rule{0mm}{5.5mm}Monte Carlo V@R&
\hs\raisebox{-1.2mm}{\includegraphics{factor.7}}&
\hs\raisebox{-1.2mm}{\includegraphics{factor.7}}&
\hs\raisebox{-1.2mm}{\includegraphics{factor.8}}&
\hs\raisebox{-1.2mm}{\includegraphics{factor.9}}&
\hs\raisebox{-1.2mm}{\includegraphics{factor.9}}&
\hs\raisebox{-1.2mm}{\includegraphics{factor.8}}&
\hs\raisebox{-1.2mm}{\includegraphics{factor.7}}&
\hs\raisebox{-1.2mm}{\includegraphics{factor.8}}&
\hs\raisebox{-1.2mm}{\includegraphics{factor.9}}\\[2mm]
\hline
\rule{0mm}{5.5mm}Historic V@R&
\hs\raisebox{-1.2mm}{\includegraphics{factor.7}}&
\hs\raisebox{-1.2mm}{\includegraphics{factor.7}}&
\hs\raisebox{-1.2mm}{\includegraphics{factor.9}}&
\hs\raisebox{-1.2mm}{\includegraphics{factor.8}}&
\hs\raisebox{-1.2mm}{\includegraphics{factor.9}}&
\hs\raisebox{-1.2mm}{\includegraphics{factor.9}}&
\hs\raisebox{-1.2mm}{\includegraphics{factor.7}}&
\hs\raisebox{-1.2mm}{\includegraphics{factor.7}}&
\hs\raisebox{-1.2mm}{\includegraphics{factor.8}}\\[2mm]
\hline
\end{tabular}

\begin{picture}(1,1)(79,-57)
\put(0,0){\includegraphics{factor.10}}
\put(2,3){\bf Technique}
\put(32,8){\rotatebox{90}{\bf Properties}}
\end{picture}

\vspace{5mm}
\parbox{130mm}{\textbf{Table~1.} Comparison of various
risk measurement techniques.
The relative strengths of different
methods are shown by the extent of shading.}
\end{center}

\skb
\textbf{4. Risk management.}
Two practical recommendations of this paper concerning
risk management are as follows:
\begin{mitemize}
\item[\bf 1.] The central management of a firm should impose the
limits not on the outstanding risks
(resp., factor risks) of the desks' portfolios,
but rather on their risk contributions
(resp., factor risk contributions)
to the capital of the whole firm.
In view of the formulas given above, this can be done
simply by announcing the arrays $(t_{kl})$ and
$(l_k)$ (resp., $(t_{kl})$ and $(l_k^m)$).
\item[\bf 2.] If the desks are allowed to trade these
risk limits within the firm, then the corresponding
competitive optimum is in fact the global optimum for
the whole firm, regardless of the initial allocation
of the risk limits.
\end{mitemize}

%======================================================
\section*{Appendix: $L^0$ Version of the Kusuoka Theorem}
\label{A}

\textbf{Theorem~A.1.}
\textit{Let $(\Omega,\F,\PP)$ be atomless and $u$ be a coherent
utility on $L^0$ with the determining set $\DDD$.
The following conditions are equivalent:
\begin{mitemize}
\item[\rm(i)] $u$ is law invariant;
\item[\rm(ii)] $\DDD$ is law invariant {\rm(}i.e.
if $Z\in\DDD$ and $Z'\stackrel{\text{\rm Law}}{=}Z$,
then $Z'\in\DDD${\rm);}
\item[\rm(iii)] there exists a set $\mathfrak M$
of probability measures on $(0,1]$ such that
$$
u(X)=\inf_{\mu\in\mathfrak M}u_\mu(X),\quad X\in L^0;
$$
\item[\rm(iv)] there exists a set $\mathfrak M$ of
probability measures on $(0,1]$ such that
$\DDD=\bigcup_{\mu\in\mathfrak M}\DDD_\mu$.
\end{mitemize}
}

\skm
\Proof
In essence, the reasoning will follow the lines of the
proof of the $L^\infty$-version of the theorem
borrowed from~\cite[Sect.~4.5]{FS04}.

(i)$\Rightarrow$(ii)
As $(\Omega,\F,\PP)$ is atomless, it supports a
random variable~$X$ with any given distribution~$\QQ$
on the real line. Then the equality $\wt u(\QQ):=u(X)$
correctly defines a function on the set of distributions.
By the definition,
\begin{align*}
\DDD
&=\{Z:Z\ge0,\;\EE Z=1,
\text{ and }\EE XZ\ge u(X)\;\forall X\in L^0\}\\
&=\bigcap_\QQ\{Z:Z\ge0,\;\EE Z=1,
\text{ and }\EE XZ\ge\wt u(\QQ)\;\forall X
\text{ with }\Law X=\QQ\},
\end{align*}
where the intersection is taken over all the probability
measures on~$\R$.
By the Hardy--Littlewood inequality
(see~\cite[Th.~A.24]{FS04}),
\begin{align*}
&\EE X^+Z
\ge\int_0^1 q_s(X^+)q_{1-s}(Z)ds
=\int_{F(0)}^1 q_s(X)q_{1-s}(Z)ds,\\
&\EE X^-Z
\le\int_0^1q_s(X^-)q_s(Z)ds
=\int_0^1 q_{1-s}(X^-)q_{1-s}(Z)ds
=-\int_0^{F(0)}q_s(X)q_{1-s}(Z)ds,
\end{align*}
where $F$ is the distribution function of~$X$.
Consequently,
$$
\EE XZ\ge\int_0^1 q_s(X)q_{1-s}(Z)ds.\eqno\text{(a.1)}
$$
Fix $Z$ and a probability measure~$\QQ$ on~$\R$.
As $(\Omega,\F,\PP)$ is atomless,
$Z$ can be represented as $q_U(Z)$ with a uniformly
distributed on $[0,1]$ random variable~$U$.
Then $X:=q_{1-U}(\QQ)$ has the distribution~$\QQ$
($q_\la(\QQ)$ is the $\la$-quantile of~$\QQ$).
Obviously,
$$
\EE XZ=\int_0^1 q_s(\QQ)q_{1-s}(Z)ds.
$$
Combining this with~(a.1), we get
$$
\inf_{X:\Law X=\QQ}\EE ZX=\int_0^1 q_s(\QQ)q_{1-s}(Z)ds.
$$
Thus, the set
$\{Z:\EE ZX\ge\wt u(\QQ)\;\forall X\text{ with }\Law X=\QQ\}$
is law invariant. Hence, $\DDD$ is law invariant.

(ii)$\Rightarrow$(iii)
The law invariance of $\DDD$ means that there exists
a set $\mathfrak Q$ of probability measures on~$\R_+$
such that $\DDD=\{Z:\Law Z\in\mathfrak Q\}$.
Fix $\QQ\in\mathfrak Q$.
There exists a unique measure~$\mu=\mu(\QQ)$ on
$(0,1]$ such that $\psi_\mu(x)=q_{1-x}(\QQ)$, $x\in(0,1]$
($\psi_\mu$ is given by~\eqref{crm7}).
Clearly, $\mu$ is positive and
\begin{align*}
\mu((0,1])
&=\int_{(0,1]}\int_{(0,y]}y^{-1}dx\mu(dy)
=\int_0^1\int_{[x,1]}y^{-1}\mu(dy)dx\\
&=\int_0^1 q_{1-x}(\QQ)dx
=\int_{\R_+}x\QQ(dx)
=1.
\end{align*}
Applying the same argument as above and recalling~\eqref{crm6},
we get
$$
\inf_{Z:\Law Z=\QQ}\EE XZ
=\int_0^1 q_s(X)q_{1-s}(\QQ)ds
=\int_0^1 q_s(X)\psi_\mu(X)ds
=u_\mu(X),\quad X\in L^0.\eqno\text{(a.2)}
$$
As a result,
$$
u(X)
=\inf_{Z\in\DDD}\EE XZ
=\inf_{\QQ\in\mathfrak Q}\inf_{Z:\Law Z=\QQ}\EE XZ
=\inf_{\QQ\in\mathfrak Q}u_{\mu(\QQ)}(X),\quad X\in L^0.
\eqno\text{(a.3)}
$$

(iii)$\Rightarrow$(i) This implication follows from
the law invariance of $u_\mu$.

(ii)$\Rightarrow$(iv)
It follows from~(a.2) that
$\{Z:\Law Z=\QQ\}\subseteq\DDD_{\mu(\QQ)}$. Hence,
$\DDD\subseteq\bigcup_{\QQ\in\mathfrak Q}\DDD_{\mu(\QQ)}$.
On the other hand, if $Z\in\DDD_{\mu(\QQ)}$ with some
$\QQ\in\mathfrak Q$, then, by~(a.3),
$$
\EE XZ
\ge u_{\mu(\QQ)}(X)
\ge u(X),\quad X\in L^0,
$$
so that, by definition, $Z\in\DDD$. Hence,
$\DDD=\bigcup_{\QQ\in\mathfrak Q}\DDD_{\mu(\QQ)}$.

(iv)$\Rightarrow$(ii) This implication follows from the
law invariance of $\DDD_\mu$, which is seen
from~\eqref{crm12}.\End

\skm
\textbf{Corollary~A.2.}
\textit{Let $(\Omega,\F,\PP)$ be atomless,
$u$ be a law invariant coherent utility
with the determining set~$\DDD$,
and $Y$ be a random vector.
Then $\EE(\DDD\cond Y)\subseteq\DDD$.
}

\skm
\Proof
By Theorem~A.1,
$\DDD=\bigcup_{\mu\in\mathfrak M}\DDD_\mu$ with some set
$\mathfrak M$ of probability measures on $(0,1]$.
Using representation~\eqref{crm12}, the convexity of the
function~$\Phi_\mu$ given by~\eqref{crm13}, and the
Jensen inequality, we get
$\EE(\DDD_\mu\cond Y)\subseteq\DDD_\mu$,
which yields the desired statement.\End

\skm
The following example shows that the law invariance
condition in the preceding corollary is essential.

\skm
\textbf{Example~A.3.}
Let $\DDD=\{Z\}$, where $Z\ne1$.
Let $Y$ be independent of $d\QQ/d\PP$. Then
$\EE(\DDD\cond Y)=\{1\}\not\subseteq\DDD$.\End

%======================================================

\end{document}